\documentclass[draft,preprint,12pt,numbers,sort&compress]{elsarticle}
\usepackage{}
\usepackage{mathrsfs}
\usepackage{amssymb}
\usepackage{amsthm}
\usepackage{mathrsfs}
\usepackage[centertags]{amsmath}
\usepackage{amsfonts}

\usepackage{color}

\input amssym.def
\input amssym.tex

\date{}

\newtheorem{Theorem}{Theorem}[section]

\newtheorem{Lemma}{Lemma}[section]

\newcommand\R{\mbox{\bf R}}

\newcommand\SR{\mbox{\scriptsize\bf R}}

\newcommand{\definition}{{\lower .5ex
  \hbox{$\>\>\stackrel{\triangle}{=}\>\>$} }}

\begin{document}

\baselineskip=22pt
\thispagestyle{empty}

\mbox{}
\bigskip

\begin{center}

{\Large \bf Global well-posedness of the  generalized KP-II}\\[1ex]
{\Large \bf   in anisotropic   Sobolev spaces}\\[1ex]

{Wei Yan\footnote{Email:yanwei19821115@sina.cn}$^{a}$,\quad Yongsheng Li\footnote{Email:yshli@scut.edu.cn}$^b$\quad, Yimin Zhang
\footnote{ Email: zhangyimin@whut.edu.cn}$^c$}\\[1ex]

{$^a$School of Mathematics and Information Science and  Henan Engineering Laboratory for
 Big Data Statistical Analysis and Optimal Control, Henan  Normal University,}\\
{Xinxiang, Henan 453007, P. R. China}\\[2ex]

{$^b$School of Mathematics, South  China  University of  Technology,}\\
{Guangzhou, Guangdong 510640,  P. R. China}\\[1ex]

\end{center}

\noindent{\bf Abstract.}In this paper,
 we consider the Cauchy problem
  for the generalized  KP-II  equation
\begin{eqnarray*}
 u_{t}-|D_{x}|^{\alpha}u_{x}+\partial_{x}^{-1}\partial_{y}^{2}u+\frac{1}{2}\partial_{x}(u^{2})=0,\alpha\geq4.
 \end{eqnarray*}
 The goal of  this paper is two-fold. Firstly,
 we prove that the  problem
  is locally well-posed  in  anisotropic Sobolev spaces $H^{s_{1},\>s_{2}}(\R^{2})$
  with $s_{1}>-\frac{3\alpha-2}{8}$, $s_{2}\geq 0$ and $\alpha\geq4$. Secondly,  we prove that
the problem is globally well-posed in  anisotropic Sobolev spaces $H^{s_{1},\>0}(\R^{2})$ with $s_{1}>-\frac{(3\alpha-4)^{2}}{28\alpha}$ and $\alpha\geq4$. Thus, our global well-posedness result improves the global well-posedness result of Hadac  (Transaction of the American Mathematical Society, 360(2008), 6555-6572.) when $4\leq \alpha\leq6.$

\bigskip
\bigskip

{\large\bf 1. Introduction}
\bigskip

\setcounter{Theorem}{0} \setcounter{Lemma}{0}

\setcounter{section}{1}

In this paper, we consider the Cauchy problem for the generalized  KP-I equation
\begin{eqnarray}
&& u_{t}-|D_{x}|^{\alpha}u_{x}+\partial_{x}^{-1}\partial_{y}^{2}u+\frac{1}{2}\partial_{x}(u^{2})=0,\label{1.01}\\
&&u(x,y,0)=u_{0}(x,y)\label{1.02}
\end{eqnarray}
in anisotropic Sobolev space $H^{s_{1},s_{2}}(\R^{2})$ defined in page 6. Here $\partial_{x}^{-1}$ is defined by its Fourier multiplier $-i\xi^{-1}.$
(\ref{1.01}) occurs in the modeling of certain long dispersive waves
 \cite{AS,KB1991,KB-1991}.
When $\alpha=2$,  (\ref{1.01}) reduces to  the KP-II equation
\begin{eqnarray}
 u_{t}+ \partial_{x}^{3}u+\partial_{x}^{-1}\partial_{y}^{2}u
 +\frac{1}{2}\partial_{x}(u^{2})=0.\label{1.03}
\end{eqnarray}
The  KP-II equations arise in physical contexts as models for the propagation
of dispersive long waves with weak transverse
effects \cite{KP}, which are two-dimensional extensions
 of the Korteweg-de-Vries
equation.

  Many people have investigated the Cauchy problem for KP-II equation,
 for instance, see \cite{Bourgain-GAFA-KP,BS,Hadac2008,Hadac2009,
 HN,IMS-NA,IMS1992,IMS-JMAA,IMS,IMCPDE,IMEJDE,IMT,IKT,ILP,ILM-CPDE,IM2011,LPS,KL,LJDE,ST-CMP,TDCDS,TADE,TT,TCPDE,TDIE,TIMRN}
 and the references therein. Bourgain \cite{Bourgain-GAFA-KP} established
   the global well-posedness of  the Cauchy problem for the
 KP-II equation  in $L^{2}(\R^{2})$ and $L^{2}(\mathbf{T}^{2}).$ Takaokao and  Tzvetkov \cite{TT} and
 Isaza and  Mej\'{\i}a \cite{IMCPDE} established the local well-posedness
  of KP-II equation in $H^{s_{1},s_{2}}(\R^{2})$
 with $s_{1}>-\frac13$ and $s_{2}\geq0.$ Takaokao \cite{TDCDS} established the local well-posedness
 of KP-II equation in $H^{s_{1},0}(\R^{2})$
 with $s_{1}>-\frac12$ under the assumption that $D_{x}^{-\frac{1}{2}+\epsilon}u_{0}\in L^{2}$
 with the suitable chosen $\epsilon$, where $D_{x}^{-\frac{1}{2}+\epsilon}$ is Fourier multiplier
 operator with multiplier $|\xi|^{-\frac{1}{2}+\epsilon}.$
 Hadac  et al. \cite{Hadac2009}  established the small data global well-posedness
 and scattering result of  KP-II equation in the homogeneous anisotropic Sobolev space
$\dot{H}^{-\frac{1}{2},\>0}(\R^{2})$ and arbitrary large initial data local well-posedness
 in both homogeneous Sobolev space $\dot{H}^{-\frac{1}{2},\>0}(\R^{2})$ and inhomogeneous anisotropic  Sobolev space
$H^{-\frac{1}{2},\>0}(\R^{2})$.

Some authors have studied the Cauchy problem for KP-I equation
\begin{eqnarray}
&& u_{t}-\partial_{x}^{3}u+\partial_{x}^{-1}\partial_{y}^{2}u+\frac{1}{2}\partial_{x}(u^{2})=0,\label{1.04}
\end{eqnarray}
for instance, see \cite{CIKS,CKS-GAFA,GPW,IKT,Kenig,MST2002,MST2004,MST-Duke, MST2007,MST2011,Z} and the references therein.
It is worth noticing that the resonant function of KP-I equation doesnot possess the
 good property as the KP-II equation.

When $\alpha=4,$ (\ref{1.01}) reduces to the fifth-order KP-II equation
\begin{eqnarray}
&& u_{t}-\partial_{x}^{5}u+\partial_{x}^{-1}\partial_{y}^{2}u+\frac{1}{2}\partial_{x}(u^{2})=0.\label{1.05}
\end{eqnarray}
Saut and Tzvetkov \cite{ST1999} proved that the Cauchy problem for (\ref{1.05}) is locally well-posed
in $H^{s_{1},s_{2}}(\R^{2})$ with $s_{1}>-\frac14,s_{2}\geq0$. Isaza et al. \cite{ILM-CPAA} proved that the
Cauchy problem for (\ref{1.05}) is locally well-posed
in $H^{s_{1},s_{2}}(\R^{2})$ with $s_{1}>-\frac54,s_{2}\geq0$ and globally well-posed
in $H^{s_{1},0}(\R^{2})$ with $s_{1}>-\frac47.$ Recently, Li and Shi \cite{LS} proved that
the Cauchy problem for
(\ref{1.05}) is locally well-posed
in $H^{s_{1},s_{2}}(\R^{2})$ with $s_{1}\geq-\frac54,s_{2}\geq0$.

Some people have studied the Cauchy problem for the fifth order KP-I equation
\begin{eqnarray}
&& u_{t}+\partial_{x}^{5}u+\partial_{x}^{-1}\partial_{y}^{2}u+\frac{1}{2}\partial_{x}(u^{2})=0,\label{1.06}
\end{eqnarray}
for instance, see \cite{ST2000,CLM,LX,ST1999,GHF} and the references therein.

In this paper, inspired by \cite{CKS-GAFA,ST2000,LX,ILM-CPAA},
by using the Fourier restriction norm method introduced in \cite{Beals,Bourgain93,KM,RR}
and developed in \cite{KPV1993,KPV1996}, the Cauchy-Schwartz inequality and Strichartz estimates
 as well as  suitable splitting of domains,  we prove
 that the Cauchy problem for (\ref{1.01})
is locally well-posed  in the anisotropic Sobolev spaces $H^{s_{1},\>s_{2}}(\R^{2})$
  with $s_{1}>-\frac{3\alpha-2}{8}$ and $s_{2}\geq 0$; using  the local well-posedness result
  of this paper and  the I-method appeared in \cite{CKSTT2001,CKSTT2003},
   we also prove that
the problem is globally well-posed in $H^{s_{1},\>0}(\R^{2})$ with $s_{1}>-\frac{(3\alpha-4)^{2}}{28\alpha}$.
Thus, our result improves the result of \cite{Hadac2008}.

We introduce some notations before presenting the main results.
 Throughout this paper, we assume that
$C$ is a positive constant which may depend upon $\alpha$
 and  vary from line to line. $a\sim b$ means that there exist constants $C_{j}>0(j=1,2)$ such that $C_{1}|b|\leq |a|\leq C_{2}|b|$.
 $a\gg b$ means that there exist a positive constant $C^{\prime}$ such that  $|a|> C^{\prime}|b|.$ $0<\epsilon\ll1$ means that $0<\epsilon<\frac{1}{100\alpha}$.
 We define
  \begin{eqnarray*}
  &&\langle\cdot\rangle:=1+|\cdot|,\\
  &&\phi(\xi,\mu):=-\xi|\xi|^{\alpha}+\frac{\mu^{2}}{\xi},\\
  &&\sigma:=\tau-\phi(\xi,\mu),\sigma_{j}=\tau_{j}-\phi(\xi_{j},\mu_{j})(j=1,2),\\
  &&\mathscr{F}u(\xi,\mu,\tau):=\frac{1}{(2\pi)^{\frac{3}{2}}}\int_{\SR^{3}}e^{-ix\xi-iy\mu-it\tau}u(x,y,t)dxdydt,\\
  &&\mathscr{F}_{xy}f(\xi,\mu):=\frac{1}{2\pi}\int_{\SR^{2}}e^{-ix\xi-iy\mu}f(x,y)dxdy,\\
  &&\mathscr{F}^{-1}u(\xi,\mu,\tau):=\frac{1}{(2\pi)^{\frac{3}{2}}}\int_{\SR^{3}}e^{ix\xi+iy\mu+it\tau}u(x,y,t)dxdydt,\\
  &&D_{x}^{a}u(x,y,t):=\frac{1}{(2\pi)^{\frac{3}{2}}}\int_{\SR^{2}}|\xi|^{a}\mathscr{F}u(\xi,\mu,\tau)e^{ix\xi+iy\mu+it\tau}d\xi d\mu d\tau,\\
 && W(t)f:=\frac{1}{2\pi}\int_{\SR^{2}}e^{ix\xi+iy\mu+it\phi(\xi,\mu)}\mathscr{F}_{xy}f(\xi,\mu)d\xi d\mu.
  \end{eqnarray*}
   Let $\eta$ be a bump function with compact support in $[-2,2]\subset \R$
  and $\eta=1$ on $(-1,1)\subset \R$.
 For each integer $j\geq1$, we define $\eta_{j}(\xi)=\eta(2^{-j}\xi)-\eta(2^{1-j}\xi),$
  $\eta_{0}(\xi)=\eta(\xi),$
 $\eta_{j}(\xi,\mu,\tau)=\eta_{j}(\sigma),$ thus, $\sum\limits_{j\geq0}\eta_{j}(\sigma)=1.$
 $\psi(t)$ is a smooth function
 supported in $[0,2]$ and equals
 $1$ in $[0,1]$.

\noindent We define
 \begin{eqnarray*}
 \|f\|_{L_{t}^{r}L_{xy}^{p}}:=\left(\int_{\SR}\left(\int_{\SR^{2}}|f|^{p}dxdy\right)^{\frac{r}{p}}dt\right)^{\frac{1}{r}}.
 \end{eqnarray*}
 For $s_{1},s_{2}\in \R,$ the anisotropic Sobolev space $H^{s_{1},s_{2}}$ is defined as follows:
 \begin{eqnarray*}
 H^{s_{1},s_{2}}(\R^{2}):=\left\{u_{0}\in \mathscr{S}^{'}(\R^{2}):\quad \|u_{0}\|_{H^{s_{1},s_{2}}(\SR^{2})}=\left\|\langle\xi\rangle^{s_{1}}
 \langle\mu\rangle^{s_{2}}\mathscr{F}_{xy}u_{0}(\xi,\mu)\right\|_{L_{\xi\mu}^{2}}<\infty\right\}
 \end{eqnarray*}
 and
 space
$
  X_{b}^{s_{1},s_{2}}
$ is defined by
$$
X_{b}^{s_{1},s_{2}}:= \left\{u\in  \mathscr{S}^{'}(\R^{3})
 :\, \|u\|_{X_{b}^{s_{1},s_{2}}}
 =  \left\|\langle\xi\rangle^{s_{1}} \langle\mu\rangle^{s_{2}}
 \left\langle\sigma\right\rangle^{b}\mathscr{F}u(\xi,\mu,\tau)
 \right\|_{L_{\tau\xi\mu}^{2}(\SR^{3})}<\infty\right\}.
$$
The space $ X_{b}^{s_{1},s_{2}}([0,T])$ denotes the restriction
 of $X_{b}^{s_{1},s_{2}}$ onto the finite time interval $[0,T]$ and
is equipped with the norm
 \begin{equation*}
    \|u\|_{X_{b}^{s_{1},s_{2}}([0,T])} =\inf \left\{\|g\|_{X_{b}^{s_{1},s_{2}}}
    :g\in X_{b}^{s_{1},s_{2}}, u(t)=g(t)
 \>\> {\rm for} \>  t\in [0,T]\right\}.
 \end{equation*}
For $s<0$ and $N\in N^{+}$, $N\geq 100$, we define an operator $I_{N}$  by
 $\mathscr{F}I_{N}u(\xi,\mu,\tau)=M(\xi)\mathscr{F}u(\xi,\mu,\tau)$,
where $M(\xi)=1$ if $|\xi|<N$; $M(\xi)=\left(\frac{|\xi|}{N}\right)^{s}$ if $|\xi|\geq N.$

The main results of this paper are as follows.

\begin{Theorem}\label{Thm1}(Local well-posedness)

\noindent Let $\alpha\geq4$ and  $|\xi|^{-1}\mathscr{F}_{xy}u_{0}(\xi,\mu)\in \mathscr{S}^{'}(\R^{2})$.  Then,  the Cauchy problem  for (\ref{1.01}) is  locally well-posed in
 $H^{s_{1},\>s_{2}}(\R^{2})$ with $s_{1}>\frac{1}{4}-\frac{3}{8}\alpha,\>s_{2}\geq0.$
\end{Theorem}
\noindent{\bf Remark 1.} When $4\leq \alpha\leq6,$ Hadac \cite{Hadac2008} has proved that the Cauchy problem for (\ref{1.01}) is locally well-posed
in  $H^{s_{1},\>s_{2}}(\R^{2})$ with $s_{1}>\frac{1}{4}-\frac{3}{8}\alpha,\>s_{2}\geq0.$ Thus, our result extends the result of Hadac \cite{Hadac2008}.

\begin{Theorem}\label{Thm2}(Global well-posedness)

\noindent
 Let $\alpha\geq4$ and  $|\xi|^{-1}\mathscr{F}_{xy}u_{0}(\xi,\mu)\in \mathscr{S}^{'}(\R^{2})$.
  Then, the Cauchy problem  for (\ref{1.01}) is  globally well-posed in $H^{s_{1},\>0}(\R^{2})$
with $s_{1}>-\frac{(3\alpha-4)^{2}}{28\alpha}$.
\end{Theorem}

\noindent{\bf Remark 2.}When $\alpha=4,$ we  have proved that the Cauchy problem for (\ref{1.01}) is globally well-posed in
$H^{s_{1},\>s_{2}}(\R^{2})$ with $s_{1}>-\frac{4}{7}$.  Isaza and  Mej\'{\i}a \cite{IM2011} have proved the same result of the Cauchy problem for (\ref{1.01})
when $\alpha=4$. When $4\leq \alpha\leq6,$ Hadac \cite{Hadac2008} has proved that the Cauchy problem for (\ref{1.01}) is globally well-posed
in  $H^{s_{1},\>s_{2}}(\R^{2})$ with $s_{1}\geq0,\>s_{2}=0.$ Thus, our result improves the result of Hadac \cite{Hadac2008}
when $4\leq \alpha \leq6.$ Since we can easily prove that the Cauchy problem  for (\ref{1.01}) is  globally well-posed in $H^{s_{1},\>0}(\R^{2})$
with $s_{1}\geq0$ with the aid of $L^{2}$ conservation law of (\ref{1.01}), thus, we only consider the case $-\frac{(3\alpha-4)^{2}}{28\alpha}<s_{1}<0.$

The rest of the paper is arranged as follows. In Section 2,  we give some
preliminaries. In Section 3, we establish two crucial bilinear estimates.
 In Section 4, we prove
the Theorem 1.1. In Section 5, we prove
the Theorem 1.2.

\bigskip
\bigskip

 \noindent{\large\bf 2. Preliminaries }

\setcounter{equation}{0}

\setcounter{Theorem}{0}

\setcounter{Lemma}{0}

\setcounter{section}{2}

In this section, motivated by \cite{Bourgain-GAFA-KP, MST2011}, we give
   Lemmas 2.1-2.6 which play a significant role in establishing
  Lemmas 3.1, 3.2. Lemma 2.2 in combination with Lemma 3.1  yields Theorem 1.1. Lemma 2.7 in combination with Lemmas 3.1, 3.2 yields
  Lemma 5.1.

\begin{Lemma}\label{Lemma2.1}
 Let $b>|a|\geq 0$. Then, we have that
 \begin{eqnarray}
 &&\int_{-b}^{b}\frac{dx}{\langle x+a\rangle^{\frac{1}{2}}}\leq
  Cb^{\frac{1}{2}},\label{2.01}\\
 &&\int_{\SR}\frac{dt}{\langle t\rangle^{\gamma}\langle t-a\rangle^{\gamma}}
 \leq C\langle a\rangle^{-\gamma},\gamma>1,\label{2.02}\\
 &&\int_{\SR}\frac{dt}{\langle t\rangle^{\gamma}|t-a|^{\frac{1}{2}}}
 \leq C\langle a\rangle^{-\frac{1}{2}},\gamma>1,\label{2.03}\\
 &&\int_{-K}^{K}\frac{dx}{|x|^{\frac{1}{2}}|a-x|^{\frac{1}{2}}}
 \leq C\frac{K^{\frac{1}{2}}}{|a|^{\frac{1}{2}}}.\label{2.04}
 \end{eqnarray}
 \end{Lemma}

  The  conclusion of (\ref{2.01}) is given in (2.4) of Lemma 2.1 in \cite{ILM-CPAA}.
(\ref{2.02})-(\ref{2.03})  can be seen in Proposition 2.2 of \cite{ST2000}. (\ref{2.04})
can be seen in line 24 of page 6562 in \cite{Hadac2008}.

\begin{Lemma}\label{Lemma2.2}
Let $T\in (0,1)$ and $s_{1},s_{2} \in \R$ and $-\frac{1}{2}<b^{\prime}\leq0\leq b\leq b^{\prime}+1$ and $\psi(t)$ be defined as in line 2 from bottom of page 5 .
Then, for $\phi\in H^{s_{1},\>s_{2}}$ and $h\in X_{b^{\prime}}^{s_{1},s_{2}},$  we have that
\begin{eqnarray}
&&\left\|\psi(t)W(t)\phi\right\|_{X_{b}^{s_{1},s_{2}}}\leq C\|\phi\|_{H^{s_{1},\>s_{2}}},\label{2.05}\\
&&\left\|\psi\left(\frac{t}{T}\right)\int_{0}^{t}W(t-\tau)h(\tau)d\tau\right\|_{X_{b}^{s_{1},\>s_{2}}}\leq C
T^{1+b^{\prime}-b}\|h\|_{X_{b^{\prime}}^{s_{1},\>s_{2}}}.\label{2.06}
\end{eqnarray}
\end{Lemma}

For the proof of Lemma 2.2, we refer the readers to \cite{G,Bourgain93,KPV1993}
 and  Lemmas 1.7, 1.9 of  \cite{Grunrock}.

\begin{Lemma}\label{Lemma2.3}
Let
$b>\frac{1}{2}$ and
$
G(\xi_{1},\mu_{1},\tau_{1},\xi,\mu,\tau)=f_{1}(\xi_{1},\mu_{1},\tau_{1})
f_{2}(\xi-\xi_{1},\mu-\mu_{1},\tau-\tau_{1})f(\xi,\mu,\tau),
 $
 then,  we have that
 \begin{eqnarray}
 \|u_{1}u_{2}\|_{L^{2}}\leq C\left(\prod\limits_{j=1}^{2}\|D_{x}^{\frac{1}{4}-\frac{\alpha}{8}}u_{j}\|_{X_{b}^{0,0}}\right)\label{2.07}
 \end{eqnarray}
and
\begin{eqnarray}
&&\left|\int_{\SR^{6}}\frac{|\xi_{1}|^{-\frac{1}{4}+\frac{\alpha}{8}}
|\xi-\xi_{1}|^{-\frac{1}{4}+\frac{\alpha}{8}}G(\xi_{1},\mu_{1},\tau_{1},\xi,\mu,\tau)}
{\prod\limits_{j=1}^{2}\langle\sigma_{j}\rangle^{b}}d\xi_{1}d\mu_{1}d\tau_{1}d\xi d\mu d\tau\right|\nonumber\\&&
\leq C\|f\|_{L_{\tau\xi\mu}^{2}}\left(\prod\limits_{j=1}^{2}\|f_{j}\|_{L_{\tau\xi\mu}^{2}}\right)\label{2.08}
 \end{eqnarray}
and
\begin{eqnarray}
&&\left|\int_{\SR^{6}}\frac{|\xi_{1}|^{-\frac{1}{4}+\frac{\alpha}{8}}
|\xi|^{-\frac{1}{4}+\frac{\alpha}{8}}G(\xi_{1},\mu_{1},\tau_{1},\xi,\mu,\tau)}
{\langle\sigma_{1}\rangle^{b}\langle\sigma\rangle^{b}}d\xi_{1}d\mu_{1}d\tau_{1}d\xi d\mu d\tau\right|\nonumber\\&&
\leq C\|f\|_{L_{\tau\xi\mu}^{2}}\left(\prod\limits_{j=1}^{2}\|f_{j}\|_{L_{\tau\xi\mu}^{2}}\right)\label{2.09}
 \end{eqnarray}
and
\begin{eqnarray}
&&\left|\int_{\SR^{6}}\frac{|\xi|^{-\frac{1}{4}+\frac{\alpha}{8}}
|\xi-\xi_{1}|^{-\frac{1}{4}+\frac{\alpha}{8}}G(\xi_{1},\mu_{1},\tau_{1},\xi,\mu,\tau)}
{\langle\sigma\rangle^{b}\langle\sigma_{2}\rangle^{b}}d\xi_{1}d\mu_{1}d\tau_{1}d\xi d\mu d\tau\right|\nonumber\\&&
\leq C\|f\|_{L_{\tau\xi\mu}^{2}}\left(\prod\limits_{j=1}^{2}\|f_{j}\|_{L_{\tau\xi\mu}^{2}}\right).\label{2.010}
\end{eqnarray}
\end{Lemma}

For the proof of Lemma 2.3, we refer the readers to Corollary 3.2 of \cite{Hadac2008}.

\begin{Lemma}\label{Lemma2.4}
Let
$b>\frac{1}{2}$ and
$
G(\xi_{1},\mu_{1},\tau_{1},\xi,\mu,\tau)=f_{1}(\xi_{1},\mu_{1},\tau_{1})
f_{2}(\xi-\xi_{1},\mu-\mu_{1},\tau-\tau_{1})f(\xi,\mu,\tau),
 $
 we have that
\begin{eqnarray}
&&\left|\int_{\SR^{6}}\frac{|\xi_{1}|^{-\frac{1}{2}}
|\xi-\xi_{1}|^{\frac{\alpha}{4}}G(\xi_{1},\mu_{1},\tau_{1},\xi,\mu,\tau)}
{\prod\limits_{j=1}^{2}\langle\sigma_{j}\rangle^{b}}d\xi_{1}d\mu_{1}d\tau_{1}d\xi d\mu d\tau\right|\nonumber\\&&
\leq C\|f\|_{L_{\tau\xi\mu}^{2}}\left(\prod\limits_{j=1}^{2}\|f_{j}\|_{L_{\tau\xi\mu}^{2}}\right)\label{2.011}
 \end{eqnarray}
and
\begin{eqnarray}
&&\left|\int_{\SR^{6}}\frac{|\xi_{1}|^{-\frac{1}{2}}
|\xi|^{\frac{\alpha}{4}}G(\xi_{1},\mu_{1},\tau_{1},\xi,\mu,\tau)}
{\langle\sigma_{1}\rangle^{b}\langle\sigma\rangle^{b}}d\xi_{1}d\mu_{1}d\tau_{1}d\xi d\mu d\tau\right|\nonumber\\&&
\leq C\|f\|_{L_{\tau\xi\mu}^{2}}\left(\prod\limits_{j=1}^{2}\|f_{j}\|_{L_{\tau\xi\mu}^{2}}\right)\label{2.012}
 \end{eqnarray}
and
\begin{eqnarray}
&&\left|\int_{\SR^{6}}\frac{|\xi|^{-\frac{1}{2}}
|\xi-\xi_{1}|^{\frac{\alpha}{4}}G(\xi_{1},\mu_{1},\tau_{1},\xi,\mu,\tau)}
{\langle\sigma\rangle^{b}\langle\sigma_{2}\rangle^{b}}d\xi_{1}d\mu_{1}d\tau_{1}d\xi d\mu d\tau\right|\nonumber\\&&
\leq C\|f\|_{L_{\tau\xi\mu}^{2}}\left(\prod\limits_{j=1}^{2}\|f_{j}\|_{L_{\tau\xi\mu}^{2}}\right)\label{2.013}
\end{eqnarray}
and
\begin{eqnarray}
&&\left|\int_{\SR^{6}}\frac{|\xi|^{-\frac{1}{2}}
|\xi_{1}|^{\frac{\alpha}{4}}G(\xi_{1},\mu_{1},\tau_{1},\xi,\mu,\tau)}
{\langle\sigma\rangle^{b}\langle\sigma_{1}\rangle^{b}}d\xi_{1}d\mu_{1}d\tau_{1}d\xi d\mu d\tau\right|\nonumber\\&&
\leq C\|f\|_{L_{\tau\xi\mu}^{2}}\left(\prod\limits_{j=1}^{2}\|f_{j}\|_{L_{\tau\xi\mu}^{2}}\right)\label{2.014}
\end{eqnarray}
and
\begin{eqnarray}
&&\left|\int_{\SR^{6}}\frac{|\xi-\xi_{1}|^{-\frac{1}{2}}
|\xi|^{\frac{\alpha}{4}}G(\xi_{1},\mu_{1},\tau_{1},\xi,\mu,\tau)}
{\langle\sigma_{2}\rangle^{b}\langle\sigma\rangle^{b}}d\xi_{1}d\mu_{1}d\tau_{1}d\xi d\mu d\tau\right|\nonumber\\&&
\leq C\|f\|_{L_{\tau\xi\mu}^{2}}\left(\prod\limits_{j=1}^{2}\|f_{j}\|_{L_{\tau\xi\mu}^{2}}\right).\label{2.015}
 \end{eqnarray}
\end{Lemma}

For the proof of Lemma 2.4, we refer the readers to Proposition 3.5 of \cite{Hadac2008}.

\begin{Lemma}\label{Lemma2.5}
Let $\phi_{\alpha}(\xi)=\xi|\xi|^{\alpha}$, $\xi=\xi_{1}+\xi_{2}$ and $\alpha\geq4$ and
\begin{eqnarray}
r_{\alpha}(\xi,\xi_{1}):=\phi_{\alpha}(\xi)-\phi_{\alpha}(\xi_{1})-\phi_{\alpha}(\xi_{2}).\label{2.016}
\end{eqnarray}
Then $r_{\alpha}(\xi,\xi_{1})\xi\xi_{1}\xi_{2}\geq0$.
\end{Lemma}
\noindent{\bf Proof.}Hadac \cite{Hadac2008} and  Gr\"unrock et al. \cite{GPS} have given Lemma 2.7, however, they do not give the proof. Now
we give the proof.

We consider the following six cases:
\begin{eqnarray}
&&(1): \xi_{1}\geq0,\xi_{2}\geq0,\xi\geq0;\nonumber\\
&&(2): \xi_{1}\leq0,\xi_{2}\leq0,\xi\leq0;\nonumber\\
&&(3): \xi_{1}\geq0,\xi_{2}\leq0,\xi\geq0;\nonumber\\
&&(4): \xi_{1}\leq0, \xi_{2}\geq0,\xi\geq0;\nonumber\\
&&(5):\xi_{1}\geq0,\xi_{2}\leq 0, \xi\leq0;\nonumber\\
&&(6): \xi_{1}\leq0,\xi_{2}\geq0, \xi\leq0.\nonumber
\end{eqnarray}
We only consider cases (1), (3), (5) due to the symmetry.

\noindent When $\xi_{1}=0$ or $\xi_{2}=0$ or $\xi=0$, (\ref{2.012})  is valid.
Thus, we can assume that $\xi_{1}\xi_{2}\xi\neq0$. Let $x:=\frac{\xi_{1}}{\xi}.$

\noindent When (1) is valid,  we have that
\begin{eqnarray}
r_{\alpha}(\xi,\xi_{1})=\xi^{\alpha+1}\left[1-x^{\alpha+1}-(1-x)^{\alpha+1}\right].\label{2.017}
\end{eqnarray}
Here $0<x<1.$
Let
\begin{eqnarray}
F(x)=1-x^{\alpha+1}-(1-x)^{\alpha+1}.\label{2.018}
\end{eqnarray}
From (\ref{2.018}),  we have that
\begin{eqnarray}
F^{\prime}(x)=(\alpha+1)\left[(1-x)^{\alpha}-x^{\alpha}\right].\label{2.019}
\end{eqnarray}
When $0<x\leq \frac{1}{2}$, we have that
\begin{eqnarray}
F^{\prime}(x)=(\alpha+1)\left[(1-x)^{\alpha}-x^{\alpha}\right]\geq0.\label{2.020}
\end{eqnarray}
From (\ref{2.020}), we have that
\begin{eqnarray}
F(x)\geq \lim\limits_{x\rightarrow0}F(x)=0.\label{2.021}
\end{eqnarray}
When $\frac{1}{2}\leq x<1$, we have that
\begin{eqnarray}
F^{\prime}(x)=(\alpha+1)\left[(1-x)^{\alpha}-x^{\alpha}\right]\leq0.\label{2.022}
\end{eqnarray}
From (\ref{2.022}), we have that
\begin{eqnarray}
F(x)\geq \lim\limits_{x\rightarrow1}F(x)=0.\label{2.023}
\end{eqnarray}
Combining (\ref{2.017}) with  (\ref{2.021}),  (\ref{2.023}),  we have that
 $r_{\alpha}(\xi,\xi_{1})\xi\xi_{1}\xi_{2}\geq0$ is valid.

By using a  similar to case (1), we can deal with the case (3),(5).

This completes the proof of Lemma 2.5.

\begin{Lemma}\label{Lemma2.6}
Let $\xi=\xi_{1}+\xi_{2},\tau=\tau_{1}+\tau_{2}$ and $|\xi_{\rm max}|={\rm max}\left\{|\xi|,|\xi_{1}|,|\xi_{2}|\right\}$ and $|\xi_{\rm min}|={\rm min}\left\{|\xi|,|\xi_{1}|,|\xi_{2}|\right\}$.  Then, we have
\begin{eqnarray}
&&3{\rm max}\left\{|\sigma|,|\sigma_{1}|,|\sigma_{2}|\right\}\geq|\sigma-\sigma_{1}-\sigma_{2}|=|\phi(\xi,\mu)-\phi(\xi_{1},\mu_{1})-\phi(\xi_{2},\mu_{2})|
\nonumber\\&&=\left|r_{\alpha}(\xi,\xi_{1})+\frac{(\mu_{1}\xi_{2}-\mu_{2}\xi_{1})^{2}}{\xi\xi_{1}\xi_{2}}\right|\geq |r_{\alpha}(\xi,\xi_{1})|\geq C|\xi_{\rm max}|^{\alpha}|\xi_{\rm min}|.\label{2.024}
\end{eqnarray}
\end{Lemma}
\noindent {\bf Proof.}Combining Lemma 2.5 with Lemma 3.4 of \cite{Hadac2008}, we have that (\ref{2.024}) is valid.Here the last but one inequality follows from the fact that
$r_{\alpha}$ and $\xi\xi_{1}\xi_{2}$  have the same sign as it is proved in Lemma 2.7.

This completes the proof of Lemma 2.6.

\begin{Lemma}\label{Lemma2.7}
Let $0<b_{1}<b_{2}<\frac12$. Then,
\begin{eqnarray}
&&\left\|\chi_{I}(\cdot)u\right\|_{X_{b_{1}}^{0,0}}\leq C\left\|u\right\|_{X_{b_{2}}^{0,0}},\label{2.025}\\
&&\left\|\chi_{I}(\cdot)u\right\|_{X_{-b_{2}}^{0,0}}\leq C\left\|u\right\|_{X_{-b_{1}}^{0,0}},\label{2.026}
\end{eqnarray}
\end{Lemma}
where
$\chi_{I}$ denotes the characteristic function of the time interval $I:=[0,1].$

For the proof of Lemma 2.7, we refer the readers to  Lemma 3.1 of \cite{IMEJDE}.

\bigskip
\bigskip

\noindent{\large\bf 3. Bilinear estimates}

\setcounter{equation}{0}

 \setcounter{Theorem}{0}

\setcounter{Lemma}{0}

 \setcounter{section}{3}
 In this section, we give the proof of  Lemmas 3.1-3.3.
 Lemma 3.1 is used to prove Theorem 1.1. Lemma 3.2 in combination with I-method yields Theorems 1.2. Lemma 3.3 is used to prove Lemma 5.1.

 \begin{Lemma}\label{Lemma3.1}
Let $0<\epsilon\leq\frac{1}{100\alpha}$, $\alpha\geq4$ and  $s_{1}\geq\frac{1}{4}-\frac{3}{8}\alpha+4\alpha \epsilon,s_{2}\geq0$ and $u_{j}\in X_{\frac{1}{2}+\epsilon}^{s_{1},s_{2}}(j=1,2)$.
Then, we have that
\begin{eqnarray}
&&\|\partial_{x}(u_{1}u_{2})\|_{X_{-\frac{1}{2}+2\epsilon}^{s_{1},s_{2}}}\leq C
\left(\prod_{j=1}^{2}\|u_{j}\|_{X_{\frac{1}{2}+\epsilon}^{s_{1},s_{2}}}\right).\label{3.01}
\end{eqnarray}
\end{Lemma}
\noindent{\bf Proof.}  To prove (\ref{3.01}),  by duality, it suffices to  prove that
\begin{eqnarray}
&&\left|\int_{\SR^{3}}\bar{u}\partial_{x}(u_{1}u_{2})dxdydt\right|\leq
C\|u\|_{X_{\frac{1}{2}-2\epsilon}^{-s_{1},-s_{2}}}\left(\prod_{j=1}^{2}
\|u_{j}\|_{X_{\frac{1}{2}+\epsilon}^{s_{1},s_{2}}}\right).\label{3.02}
\end{eqnarray}
for $u\in X_{\frac{1}{2}-2\epsilon}^{-s_{1},-s_{2}}.$
Let
\begin{eqnarray}
&&\xi=\xi_{1}+\xi_{2},\mu=\mu_{1}+\mu_{2},\tau=\tau_{1}+\tau_{2},\nonumber\\
&&F(\xi,\mu,\tau)=\langle\xi\rangle^{-s_{1}}\langle\mu\rangle^{-s_{2}}
\langle \sigma\rangle^{\frac{1}{2}-2\epsilon}\mathscr{F}u(\xi,\mu,\tau),\nonumber\\&&
F_{j}(\xi_{j},\mu_{j},\tau_{j})=\langle\xi_{j}\rangle^{s_{1}}\langle\mu_{j}\rangle^{s_{2}}
\langle \sigma_{j}\rangle^{\frac{1}{2}+\epsilon}
\mathscr{F}u_{j}(\xi_{j},\mu,\tau_{j})(j=1,2).\label{3.03}
\end{eqnarray}
To obtain (\ref{3.02}), from (\ref{3.03}), it suffices to prove that
\begin{eqnarray}
&&\int_{\SR^{6}}\frac{|\xi|\langle\xi\rangle^{s_{1}}\langle\mu\rangle^{s_{2}}
F(\xi,\mu,\tau)\prod\limits_{j=1}^{2}F_{j}(\xi_{j},\mu_{j},\tau_{j})}{\langle\sigma\rangle^{\frac{1}{2}-2\epsilon}
\prod\limits_{j=1}^{2}\langle\xi_{j}\rangle^{s_{1}}\langle\mu_{j}\rangle^{s_{2}}\langle\sigma_{j}\rangle^{\frac{1}{2}+\epsilon}}
d\xi_{1}d\mu_{1}d\tau_{1}d\xi d\mu d\tau\nonumber\\&&\leq C
\|F\|_{L_{\tau\xi\mu}^{2}}\left(\prod_{j=1}^{2}\|F_{j}\|_{L_{\tau\xi\mu}^{2}}\right).\label{3.04}
\end{eqnarray}
Without loss of generality, by using the symmetry,  we assume that
$|\xi_{1}|\geq |\xi_{2}|$  and   $F(\xi,\mu,\tau)\geq 0,F_j(\xi_{j},\mu_{j},\tau_{j})\geq 0(j=1,2)$
and
\begin{eqnarray*}
D^{*}:=\left\{(\xi_1,\mu_{1},\tau_1,\xi,\mu,\tau)\in \R^{6},|\xi_{1}|\geq|\xi_{2}|\right\}.
\end{eqnarray*}
We define
\begin{eqnarray*}
&&\hspace{-0.8cm}\Omega_1=\left\{(\xi_1,\mu_{1},\tau_1,\xi,\mu,\tau)\in D^{*},
 |\xi_2|\leq |\xi_{1}|\leq 80\right\},\\
&&\hspace{-0.8cm} \Omega_2=\{ (\xi_1,\mu_{1},\tau_1,\xi,\mu,\tau)\in D^{*},
|\xi_1|\geq 80, |\xi_{1}|\gg|\xi_{2}|,|\xi_{2}|\leq 20\},\\
&&\hspace{-0.8cm} \Omega_3=\{ (\xi_1,\mu_{1},\tau_1,\xi,\mu,\tau)\in D^{*},
|\xi_1|\geq 80, |\xi_{1}|\gg|\xi_{2}|,|\xi_{2}|> 20\},\\
&&\hspace{-0.8cm}\Omega_4=\{(\xi_1,\mu_{1},\tau_1,\xi,\mu,\tau)\in D^{*},
|\xi_{1}|\geq 80,4|\xi|\leq |\xi_{2}|\sim|\xi_{1}|,|\xi|\leq 20,\xi_{1}\xi_{2}<0\},\\
&&\hspace{-0.8cm}\Omega_5=\{(\xi_1,\mu_{1},\tau_1,\xi,\mu,\tau)\in D^{*},
|\xi_{1}|\geq 80,4|\xi|\leq |\xi_{2}|\sim|\xi_{1}|,|\xi|> 20,\xi_{1}\xi_{2}<0\},\\
&&\hspace{-0.8cm}\Omega_6=\{(\xi_1,\mu_{1},\tau_1,\xi,\mu,\tau)\in D^{*},
 |\xi_{1}|\geq 80, |\xi_{1}|\sim |\xi_{2}|,\xi_{1}\xi_{2}<0,|\xi|\geq \frac{|\xi_{2}|}{4}\},\\
 &&\hspace{-0.8cm}\Omega_7=\left\{(\xi_1,\mu_{1},\tau_1,\xi,\mu,\tau)\in D^{*},
 |\xi_{1}|\geq 80, |\xi_{1}|\sim |\xi_{2}|,\xi_{1}\xi_{2}>0\right\}.
\end{eqnarray*}
Obviously, $D^{*}\subset\bigcup\limits_{j=1}^{7}\Omega_{j}.$
We define
\begin{equation}
    K_{1}(\xi_{1},\mu_{1},\tau_{1},\xi,\mu,\tau):=\frac{|\xi|
    \langle\xi\rangle^{s_{1}}\langle\mu\rangle^{s_{2}}}{\langle\sigma_{j}\rangle^{\frac{1}{2}-2\epsilon}
\prod\limits_{j=1}^{2}\langle\xi_{j}\rangle^{s_{1}}
\langle\mu_{j}\rangle^{s_{2}}\langle\sigma_{j}\rangle^{\frac{1}{2}+\epsilon}}\label{3.05}
\end{equation}
and
\begin{eqnarray*}
{\rm Int_{j}}:=\int_{\Omega_{j}} K_{1}(\xi_{1},\mu_{1},\tau_{1},\xi,\mu,\tau)F(\xi,\mu,\tau)
\prod_{j=1}^{2}F_{j}(\xi_{j},\mu_{j},\tau_{j})
d\xi_{1}d\mu_{1}d\tau_{1}d\xi d\mu d\tau
\end{eqnarray*}
with $1\leq j\leq 7, j\in N.$
Since $s_{2}\geq0$ and $\mu=\sum\limits_{j=1}^{2}\mu_{j}$, we have that
 $\langle \mu\rangle^{s_{2}} \leq \prod\limits_{j=1}^{2}\langle\mu_{j}\rangle^{s_{2}}$,
thus, we have that
\begin{eqnarray}
K_{1}(\xi_{1},\mu_{1},\tau_{1},\xi,\mu,\tau)\leq\frac{|\xi|\langle\xi\rangle^{s_{1}}}
{\langle\sigma\rangle^{\frac{1}{2}-2\epsilon}
\prod\limits_{j=1}^{2}\langle\xi_{j}\rangle^{s_{1}}\langle\sigma_{j}\rangle^{\frac{1}{2}+\epsilon}}.\label{3.06}
\end{eqnarray}
(1). Region $\Omega_{1}.$ In this region $|\xi|\leq |\xi_{1}|+|\xi_{2}|\leq 160,$ thus, we have that
\begin{eqnarray*}
K_{1}(\xi_{1},\mu_{1},\tau_{1},\xi,\mu,\tau)\leq\frac{|\xi|\langle\xi\rangle^{s_{1}}}
{\langle\sigma\rangle^{\frac{1}{2}-2\epsilon}
\prod\limits_{j=1}^{2}\langle\xi_{j}\rangle^{s_{1}}\langle\sigma_{j}\rangle^{\frac{1}{2}+\epsilon}}\leq \frac{C|\xi|}{\langle\sigma\rangle^{\frac{1}{2}-2\epsilon}
\prod\limits_{j=1}^{2}\langle\sigma_{j}\rangle^{\frac{1}{2}+\epsilon}},
\end{eqnarray*}
 this case can be proved similarly to case
$low+low\longrightarrow low$ of pages 344-345 of Theorem 3.1 in \cite{LX}.

 \noindent (2). Region $\Omega_{2}.$ In this region, we have that $|\xi|\sim |\xi_{1}|$. Thus, we have that
 \begin{eqnarray}
K_{1}(\xi_{1},\mu_{1},\tau_{1},\xi,\mu,\tau)\leq C\frac{|\xi|}
{\langle\sigma\rangle^{\frac{1}{2}-2\epsilon}
\prod\limits_{j=1}^{2}\langle\sigma_{j}\rangle^{\frac{1}{2}+\epsilon}}.\label{3.07}
\end{eqnarray}
By using the Cauchy-Schwartz inequality with respect to $\xi_{1},\mu_{1},\tau_{1}$,
 from (\ref{3.07}),  we have that
\begin{eqnarray}
&&{\rm Int_{2}}\leq C\int_{\SR^{3}}\frac{|\xi|}{\langle \sigma \rangle^{\frac{1}{2}-2\epsilon}}
\left(\int_{\SR^{3}}\frac{d\xi_{1}d\mu_{1}d\tau_{1}}{\prod\limits_{j=1}^{2}
\langle\sigma_{j}\rangle^{1+2\epsilon}}\right)^{\frac{1}{2}}\nonumber\\&&\qquad\qquad\qquad \times\left(\int_{\SR^{3}}
\prod\limits_{j=1}^{2}\left|F_{j}(\xi_{j},\mu_{j},\tau_{j})\right|^{2}
d\xi_{1}d\mu_{1}d\tau_{1}\right)^{\frac{1}{2}}F(\xi,\mu,\tau)d\xi d\mu d\tau.\label{3.08}
\end{eqnarray}
By using (\ref{2.02}), we have that
\begin{eqnarray}
&&\hspace{-1.4cm}\frac{|\xi|}{\langle \sigma \rangle^{\frac{1}{2}-2\epsilon}}\left(\int_{\SR^{3}}
\frac{d\xi_{1}d\mu_{1}d\tau_{1}}{\prod\limits_{j=1}^{2}
\langle\sigma_{j}\rangle^{1+2\epsilon}}\right)^{\frac{1}{2}}
\leq C\frac{|\xi|}{\langle \sigma \rangle^{\frac{1}{2}-2\epsilon}}\left(\int_{\SR^{2}}\frac{d\xi_{1}d\mu_{1}}{
\langle\tau-\phi(\xi_{1},\mu_{1})-\phi(\xi_{2},\mu_{2})\rangle^{1+2\epsilon}}\right)^{\frac{1}{2}}\label{3.09}.
\end{eqnarray}
Let $\nu=\tau-\phi(\xi_{1},\mu_{1})-\phi(\xi_{2},\mu_{2})$, since $|\xi_{1}|\gg|\xi_{2}|$,
then we have that the absolute value of Jacobian determinant equals
\begin{eqnarray}
&&\hspace{-1cm}\left|\frac{\partial(r_{\alpha}(\xi,\xi_{1}),\nu)}{\partial(\xi_{1},\mu_{1})}\right|
=2(\alpha+1)\left|\frac{\mu_{1}}{\xi_{1}}-\frac{\mu_{2}}{\xi_{2}}\right|
\left||\xi_{1}|^{\alpha}-|\xi_{2}|^{\alpha}\right|\nonumber\\
&&\hspace{-1cm}=2(\alpha+1)\left|\sigma-\nu+r_{\alpha}(\xi,\xi_{1})\right|^{\frac{1}{2}}\left|\frac{\xi}{\xi_{1}\xi_{2}}\right|^{\frac{1}{2}}
\left||\xi_{1}|^{\alpha}-|\xi_{2}|^{\alpha}\right|\nonumber\\&&\sim
\left|\sigma-\nu+r_{\alpha}(\xi,\xi_{1})\right|^{\frac{1}{2}}\left|\frac{\xi}{\xi_{1}\xi_{2}}\right|^{\frac{1}{2}}|\xi_{1}|^{\alpha}
.\label{3.010}
\end{eqnarray}
Inserting (\ref{3.010}) into (\ref{3.09}), by using (\ref{2.03}) and Lemma 3.4 of \cite{Hadac2008}, we have that
\begin{eqnarray}
&&\frac{|\xi|}{\langle \sigma \rangle^{\frac{1}{2}-2\epsilon}}\left(\int_{\SR^{3}}
\frac{d\xi_{1}d\mu_{1}d\tau_{1}}{\prod\limits_{j=1}^{2}
\langle\sigma_{j}\rangle^{1+2\epsilon}}\right)^{\frac{1}{2}}
\leq C\frac{|\xi|}{\langle \sigma \rangle^{\frac{1}{2}-2\epsilon}}\left(\int_{\SR^{2}}\frac{d\xi_{1}d\mu_{1}}{
\langle\tau-\phi(\xi_{1},\mu_{1})-\phi(\xi_{2},\mu_{2})\rangle^{1+2\epsilon}}\right)^{\frac{1}{2}}\nonumber\\
&&\leq \frac{C}{|\xi|^{\frac{\alpha}{2}-1}\langle \sigma \rangle^{\frac{1}{2}-2\epsilon}}
\left(\int_{\SR^{2}}\frac{d\nu dr_{\alpha}(\xi,\xi_{1})}{\left|\sigma-\nu+r_{\alpha}(\xi,\xi_{1})\right|^{\frac{1}{2}}
\langle\nu\rangle^{1+2\epsilon}}\right)^{\frac{1}{2}}\nonumber\\
&&\leq  \frac{C}{|\xi|^{\frac{\alpha}{2}-1}\langle \sigma \rangle^{\frac{1}{2}-2\epsilon}}
\left(\int_{|r_{\alpha}(\xi,\xi_{1})|<20\alpha|\xi|^{\alpha}}\frac{dr_{\alpha}(\xi,\xi_{1})}{
\langle r_{\alpha}(\xi,\xi_{1})+\sigma\rangle^{\frac{1}{2}}}\right)^{\frac{1}{2}}
\label{3.011}.
\end{eqnarray}
When $|\sigma|<20\alpha|\xi|^{\alpha},$ combining (\ref{3.011}) with (\ref{2.01}), since $\alpha\geq4,$ we have that
\begin{eqnarray}
 \frac{C}{|\xi|^{\frac{\alpha}{2}-1}\langle \sigma \rangle^{\frac{1}{2}-2\epsilon}}\left(\int_{|r_{\alpha}(\xi,\xi_{1})|<20\alpha|\xi|^{\alpha}}\frac{dr_{\alpha}(\xi,\xi_{1})}{
\langle r_{\alpha}(\xi,\xi_{1})+\sigma\rangle^{\frac{1}{2}}}\right)^{\frac{1}{2}}
\leq\frac{C}{|\xi|^{\frac{\alpha}{4}-1}\langle \sigma \rangle^{\frac{1}{2}-2\epsilon}}\leq C\label{3.012}.
\end{eqnarray}
When $|\sigma|\geq20\alpha|\xi|^{\alpha},$ from   (\ref{3.011}), we have that
\begin{eqnarray}
 &&\frac{C}{|\xi|^{\frac{\alpha}{2}-1}\langle \sigma \rangle^{\frac{1}{2}-2\epsilon}}\left(\int_{|r_{\alpha}(\xi,\xi_{1})|<20\alpha|\xi|^{\alpha}}\frac{dr_{\alpha}(\xi,\xi_{1})}{
\langle r_{\alpha}(\xi,\xi_{1})+\sigma\rangle^{\frac{1}{2}}}\right)^{\frac{1}{2}}\leq\frac{C|\xi|^{\frac{\alpha}{2}}}
{|\xi|^{\frac{\alpha}{2}-1}\langle \sigma \rangle^{\frac{1}{2}-2\epsilon}}\nonumber\\&&\leq\frac{C|\xi|}
{\langle \sigma \rangle^{\frac{1}{2}-2\epsilon}}\leq C\label{3.013}.
\end{eqnarray}
Combining (\ref{3.09}) with (\ref{3.010})-(\ref{3.013}), we have that
\begin{eqnarray}
&&\frac{|\xi|}{\langle \sigma \rangle^{\frac{1}{2}-2\epsilon}}\left(\int_{\SR^{3}}
\frac{d\xi_{1}d\mu_{1}d\tau_{1}}{\prod\limits_{j=1}^{2}
\langle\sigma_{j}\rangle^{1+2\epsilon}}\right)^{\frac{1}{2}}\leq C\label{3.014}.
\end{eqnarray}
Inserting (\ref{3.014}) into (\ref{3.08}), by using the Cauchy-Schwartz
 inequality with respect to $\xi,\mu,\tau$, we have that
\begin{eqnarray}
&&{\rm Int_{2}}\leq C\int_{\SR^{3}}\frac{|\xi|}{\langle \sigma \rangle^{\frac{1}{2}-2\epsilon}}
\left(\int_{\SR^{3}}\frac{d\xi_{1}d\mu_{1}d\tau_{1}}{\prod\limits_{j=1}^{2}
\langle\sigma_{j}\rangle^{1+2\epsilon}}\right)^{\frac{1}{2}}\nonumber\\&&\qquad\qquad\qquad \times\left(\int_{\SR^{3}}
\prod\limits_{j=1}^{2}\left|F_{j}(\xi_{j},\mu_{j},\tau_{j})\right|^{2}
d\xi_{1}d\mu_{1}d\tau_{1}\right)^{\frac{1}{2}}F(\xi,\mu,\tau)d\xi d\mu d\tau\nonumber\\
&&\leq C\int_{\SR^{3}}\left(\int_{\SR^{3}}
\prod\limits_{j=1}^{2}\left|F_{j}(\xi_{j},\mu_{j},\tau_{j})\right|^{2}
d\xi_{1}d\mu_{1}d\tau_{1}\right)^{\frac{1}{2}}F(\xi,\mu,\tau)d\xi d\mu d\tau\nonumber\\&&
\leq C\|F\|_{L_{\tau\xi\mu}^{2}}\left(\prod\limits_{j=1}^{2}\|F_{j}\|_{L_{\tau\xi\mu}^{2}}\right).\label{3.015}
\end{eqnarray}
\noindent (3). Region $\Omega_{3}.$ In this region, we have that $|\xi|\sim |\xi_{1}|\sim\langle \xi\rangle\sim \langle\xi_{1}\rangle$.

\noindent
Since (\ref{2.024}) is valid, we have that one of the following three cases must occur:
\begin{eqnarray}
&&|\sigma|:={\rm max}\left\{|\sigma|,|\sigma_{1}|,|\sigma_{2}|\right\}
\geq C|\xi_{1}|^{\alpha}|\xi_{2}|,\label{3.016}\\
&&|\sigma_{1}|:={\rm max}\left\{|\sigma|,|\sigma_{1}|,|\sigma_{2}|\right\}
\geq C|\xi_{1}|^{\alpha}|\xi_{2}|,\label{3.017}\\
&&|\sigma_{2}|:={\rm max}\left\{|\sigma|,|\sigma_{1}|,|\sigma_{2}|\right\}
\geq C|\xi_{1}|^{\alpha}|\xi_{2}|.\label{3.018}
\end{eqnarray}
When (\ref{3.016}) is valid, since  $s_{1}\geq\frac{1}{4}-\frac{3\alpha}{8}+4\alpha\epsilon$ and $\alpha\geq4$, we have that
\begin{eqnarray}
&&K_{1}(\xi_{1},\mu_{1},\tau_{1},\xi,\mu,\tau)\leq\frac{|\xi|\langle\xi\rangle^{s_{1}}}
{\langle\sigma\rangle^{\frac{1}{2}-2\epsilon}
\prod\limits_{j=1}^{2}\langle\xi_{j}\rangle^{s_{1}}\langle\sigma_{j}\rangle^{\frac{1}{2}+\epsilon}}\leq C\frac{|\xi|^{1-\frac{\alpha}{2}+2\alpha\epsilon}|\xi_{2}|^{-s_{1}-\frac{1}{2}+2\epsilon}}
{\prod\limits_{j=1}^{2}\langle\sigma_{j}\rangle^{\frac{1}{2}+\epsilon}}\nonumber\\&&
\leq C\frac{|\xi|^{1-\frac{\alpha}{2}+2\alpha\epsilon}|\xi_{2}|^{\frac{3\alpha}{8}-\frac{5}{8}-(4\alpha-2)\epsilon}}
{\prod\limits_{j=1}^{2}\langle\sigma_{j}\rangle^{\frac{1}{2}+\epsilon}}
\leq C\frac{|\xi_{1}|^{-\frac{1}{4}+\frac{\alpha}{8}}|\xi_{2}|^{-\frac{1}{4}+\frac{\alpha}{8}}}
{\prod\limits_{j=1}^{2}\langle\sigma_{j}\rangle^{\frac{1}{2}+\epsilon}}
.\label{3.019}
\end{eqnarray}
Thus, combining (\ref{2.08}) with (\ref{3.019}),  we have that
\begin{eqnarray*}
&&|{\rm Int_{3}}|\leq C\left|\int_{\SR^{6}}\frac{F(\xi,\mu,\tau)\prod\limits_{j=1}^{2}|\xi_{j}|^{-\frac{1}{4}+\frac{\alpha}{8}}F_{j}(\xi_{j},\mu_{j},\tau_{j})}
{\prod\limits_{j=1}^{2}\langle\sigma_{j}\rangle^{\frac{1}{2}+\epsilon}}d\xi_{1}d\mu_{1}d\tau_{1}d\xi d\mu d\tau\right|\nonumber\\&&
\leq \|F\|_{L_{\tau\xi\mu}^{2}}\left(\prod\limits_{j=1}^{2}\|F_{j}\|_{L_{\tau\xi\mu}^{2}}\right).
\end{eqnarray*}
When (\ref{3.017}) is valid, since  $s_{1}\geq\frac{1}{4}-\frac{3\alpha}{8}+4\alpha\epsilon$
 and $\langle \sigma\rangle^{-\frac{1}{2}+2\epsilon}
 \langle \sigma_{1}\rangle^{-\frac{1}{2}-\epsilon}\leq \langle \sigma\rangle^{-\frac{1}{2}-\epsilon}\langle
  \sigma_{1}\rangle^{-\frac{1}{2}+2\epsilon},$
we have that
\begin{eqnarray}
&&K_{1}(\xi_{1},\mu_{1},\tau_{1},\xi,\mu,\tau)\leq C\frac{|\xi|\langle\xi_{2}\rangle^{-s_{1}}}
{\langle\sigma_{1}\rangle^{\frac{1}{2}-2\epsilon}
\langle\sigma_{2}\rangle^{\frac{1}{2}+\epsilon}\langle\sigma\rangle^{\frac{1}{2}+\epsilon}}\leq C\frac{|\xi|^{1-\frac{\alpha}{2}+2\alpha\epsilon}|\xi_{2}|^{-s_{1}-\frac{1}{2}+2\epsilon}}
{\langle\sigma_{2}\rangle^{\frac{1}{2}+\epsilon}\langle\sigma\rangle^{\frac{1}{2}+\epsilon}}\nonumber\\&&
\nonumber\\
&&\leq C\frac{|\xi|^{1-\frac{\alpha}{2}+2\alpha\epsilon}|\xi_{2}|^{\frac{3\alpha}{8}-\frac{5}{8}-(4\alpha-2)\epsilon}}
{\langle\sigma_{2}\rangle^{\frac{1}{2}+\epsilon}\langle\sigma\rangle^{\frac{1}{2}+\epsilon}}
\leq C\frac{|\xi_{2}|^{-\frac{1}{4}+\frac{\alpha}{8}}|\xi|^{-\frac{1}{4}+\frac{\alpha}{8}}}
{\langle\sigma_{2}\rangle^{\frac{1}{2}+\epsilon}\langle\sigma\rangle^{\frac{1}{2}+\epsilon}}
.\label{3.020}
\end{eqnarray}
Thus, combining (\ref{2.010}) with (\ref{3.020}),  we have that
\begin{eqnarray*}
|{\rm Int_{3}}|\leq C\|F\|_{L_{\tau\xi\mu}^{2}}\left(\prod\limits_{j=1}^{2}\|F_{j}\|_{L_{\tau\xi\mu}^{2}}\right).
\end{eqnarray*}
When (\ref{3.018}) is valid, since  $s_{1}\geq\frac{1}{4}-\frac{3\alpha}{8}+4\alpha\epsilon$, $\alpha\geq4$,
 and $\langle \sigma\rangle^{-\frac{1}{2}+2\epsilon}\langle \sigma_{2}
 \rangle^{-\frac{1}{2}-\epsilon}\leq \langle \sigma\rangle^{-\frac{1}{2}-\epsilon}
 \langle \sigma_{2}\rangle^{-\frac{1}{2}+2\epsilon},$
we have that
\begin{eqnarray}
&&K_{1}(\xi_{1},\mu_{1},\tau_{1},\xi,\mu,\tau)\leq\frac{|\xi|\langle\xi_{2}\rangle^{-s_{1}}}
{\langle\sigma_{2}\rangle^{\frac{1}{2}-2\epsilon}\langle\sigma_{1}\rangle^{\frac{1}{2}+\epsilon}\langle\sigma\rangle^{\frac{1}{2}+\epsilon}
}\leq C\frac{|\xi|^{1-\frac{\alpha}{2}+2\alpha\epsilon}|\xi_{2}|^{-s_{1}-\frac{1}{2}+2\epsilon}}
{\langle\sigma_{1}\rangle^{\frac{1}{2}+\epsilon}\langle\sigma\rangle^{\frac{1}{2}+\epsilon}}\nonumber\\&&
\leq C\frac{|\xi|^{1-\frac{\alpha}{2}+2\alpha\epsilon}|\xi_{2}|^{\frac{3\alpha}{8}-\frac{5}{8}-(4\alpha-2)\epsilon}}
{\langle\sigma_{1}\rangle^{\frac{1}{2}+\epsilon}\langle\sigma\rangle^{\frac{1}{2}+\epsilon}}
\leq C\frac{|\xi_{1}|^{-\frac{1}{4}+\frac{\alpha}{8}}|\xi|^{-\frac{1}{4}+\frac{\alpha}{8}}}
{\langle\sigma_{1}\rangle^{\frac{1}{2}+\epsilon}\langle\sigma\rangle^{\frac{1}{2}+\epsilon}}
.\label{3.021}
\end{eqnarray}
Thus, combining (\ref{2.09}) with (\ref{3.021}),  we have that
\begin{eqnarray*}
|{\rm Int_{3}}|\leq C\|F\|_{L_{\tau\xi\mu}^{2}}\left(\prod\limits_{j=1}^{2}\|F_{j}\|_{L_{\tau\xi\mu}^{2}}\right).
\end{eqnarray*}
(4). Region $\Omega_{4}.$
In this case, we consider (\ref{3.016})-(\ref{3.018}), respectively.

\noindent When (\ref{3.016}) is valid, since  $s_{1}\geq\frac{1}{4}-\frac{3\alpha}{8}+4\alpha\epsilon$ and $\alpha\geq4$ and $|\xi|\leq 20$, we have that
\begin{eqnarray}
&&K_{1}(\xi_{1},\mu_{1},\tau_{1},\xi,\mu,\tau)\leq\frac{|\xi|\langle\xi\rangle^{s_{1}}}
{\langle\sigma\rangle^{\frac{1}{2}-2\epsilon}
\prod\limits_{j=1}^{2}\langle\xi_{j}\rangle^{s_{1}}\langle\sigma_{j}\rangle^{\frac{1}{2}+\epsilon}}\nonumber\\
&&\leq C\frac{|\xi|^{\frac{1}{2}+2\epsilon}|\xi_{2}|^{-2s_{1}-\frac{\alpha}{2}+2\alpha\epsilon}}
{\prod\limits_{j=1}^{2}\langle\sigma_{j}\rangle^{\frac{1}{2}+\epsilon}}
\leq C\frac{\prod\limits_{j=1}^{2}|\xi_{j}|^{-\frac{1}{4}+\frac{\alpha}{8}}}
{\prod\limits_{j=1}^{2}\langle\sigma_{j}\rangle^{\frac{1}{2}+\epsilon}}
.\label{3.022}
\end{eqnarray}
Thus, combining (\ref{2.08}) with (\ref{3.022}),  we have that
\begin{eqnarray*}
|{\rm Int_{4}}|\leq C\|F\|_{L_{\tau\xi\mu}^{2}}\left(\prod\limits_{j=1}^{2}\|F_{j}\|_{L_{\tau\xi\mu}^{2}}\right).
\end{eqnarray*}
When (\ref{3.017}) is valid, since  $s_{1}\geq\frac{1}{4}-\frac{3\alpha}{8}+4\alpha\epsilon$ and $\alpha\geq4$ and $|\xi|\leq 20$,
 $\langle \sigma\rangle^{-\frac{1}{2}+2\epsilon}\langle \sigma_{1}\rangle^{-\frac{1}{2}-\epsilon}
 \leq \langle \sigma\rangle^{-\frac{1}{2}-\epsilon}\langle \sigma_{1}\rangle^{-\frac{1}{2}+2\epsilon},$
we have that
\begin{eqnarray}
&&K_{1}(\xi_{1},\mu_{1},\tau_{1},\xi,\mu,\tau)\leq\frac{|\xi|\langle\xi\rangle^{s_{1}}}
{\langle\sigma\rangle^{\frac{1}{2}-2\epsilon}
\prod\limits_{j=1}^{2}\langle\xi_{j}\rangle^{s_{1}}\langle\sigma_{j}\rangle^{\frac{1}{2}+\epsilon}}\nonumber\\
&&\leq C\frac{|\xi|^{\frac{1}{2}+2\epsilon}|\xi_{2}|^{-2s_{1}-\frac{\alpha}{2}+2\alpha\epsilon}}
{\langle\sigma\rangle^{\frac{1}{2}+\epsilon}\langle\sigma_{2}\rangle^{\frac{1}{2}+\epsilon}}
\leq C\frac{|\xi|^{-\frac{1}{2}}|\xi_{2}|^{\frac{\alpha}{4}}}
{\langle\sigma\rangle^{\frac{1}{2}+\epsilon}\langle\sigma_{2}\rangle^{\frac{1}{2}+\epsilon}}
.\label{3.023}
\end{eqnarray}
Thus, combining (\ref{2.013}) with (\ref{3.023}),  we have that
\begin{eqnarray*}
|{\rm Int_{4}}|\leq C\|F\|_{L_{\tau\xi\mu}^{2}}\left(\prod\limits_{j=1}^{2}\|F_{j}\|_{L_{\tau\xi\mu}^{2}}\right).
\end{eqnarray*}
When (\ref{3.018}) is valid, this case can be proved similarly to case (\ref{3.017}) of  Region 4 with the aid of (\ref{2.015}).

\noindent (5). Region $\Omega_{5}.$
We firstly deal with $\frac{1}{4}-\frac{3\alpha}{8}+4\alpha\epsilon\leq s_{1}\leq0.$

\noindent
When  $\frac{1}{4}-\frac{3\alpha}{8}+4\alpha\epsilon\leq s_{1}\leq0$, we consider (\ref{3.016})-(\ref{3.018}), respectively.

\noindent When (\ref{3.016}) is valid,  we have that
\begin{eqnarray}
K_{1}(\xi_{1},\mu_{1},\tau_{1},\xi,\mu,\tau)\leq\frac{|\xi|\langle\xi\rangle^{s_{1}}}
{\langle\sigma\rangle^{\frac{1}{2}-2\epsilon}
\prod\limits_{j=1}^{2}\langle\xi_{j}\rangle^{s_{1}}\langle\sigma_{j}\rangle^{\frac{1}{2}+\epsilon}}
\leq C\frac{|\xi|^{s_{1}+\frac{1}{2}+2\epsilon}|\xi_{2}|^{-2s_{1}-\frac{\alpha}{2}+2\alpha\epsilon}}
{\prod\limits_{j=1}^{2}\langle\sigma_{j}\rangle^{\frac{1}{2}+\epsilon}};\label{3.024}
\end{eqnarray}
when $s_{1}+\frac{1}{2}+2\epsilon\leq0,$ from  (\ref{3.024}), since  $s_{1}\geq\frac{1}{4}-\frac{3\alpha}{8}+4\alpha\epsilon$, we have that
\begin{eqnarray}
K_{1}(\xi_{1},\mu_{1},\tau_{1},\xi,\mu,\tau)\leq C\frac{\prod\limits_{j=1}^{2}|\xi_{j}|^{-\frac{1}{4}+\frac{\alpha}{8}}}
{\prod\limits_{j=1}^{2}\langle\sigma_{j}\rangle^{\frac{1}{2}+\epsilon}};\label{3.025}
\end{eqnarray}
when $s_{1}+\frac{1}{2}+2\epsilon>0,$ from  (\ref{3.024}), since  $s_{1}\geq\frac{1}{4}-\frac{3\alpha}{8}+4\alpha\epsilon$, we have that
\begin{eqnarray}
K_{1}(\xi_{1},\mu_{1},\tau_{1},\xi,\mu,\tau)\leq C\frac{|\xi_{1}|^{-s_{1}+\frac{1-\alpha}{2}+(2\alpha+2)\epsilon}}
{\prod\limits_{j=1}^{2}\langle\sigma_{j}\rangle^{\frac{1}{2}+\epsilon}}\leq C\frac{\prod\limits_{j=1}^{2}|\xi_{j}|^{-\frac{1}{4}+\frac{\alpha}{8}}}
{\prod\limits_{j=1}^{2}\langle\sigma_{j}\rangle^{\frac{1}{2}+\epsilon}};\label{3.026}
\end{eqnarray}
Thus, combining (\ref{2.08}) with (\ref{3.025})-(\ref{3.026}),  we have that
\begin{eqnarray*}
|{\rm Int_{5}}|\leq C\|F\|_{L_{\tau\xi\mu}^{2}}\left(\prod\limits_{j=1}^{2}\|F_{j}\|_{L_{\tau\xi\mu}^{2}}\right).
\end{eqnarray*}
When (\ref{3.017}) is valid, since
 $\langle \sigma\rangle^{-\frac{1}{2}+2\epsilon}\langle \sigma_{1}\rangle^{-\frac{1}{2}-\epsilon}
 \leq \langle \sigma\rangle^{-\frac{1}{2}-\epsilon}\langle \sigma_{1}\rangle^{-\frac{1}{2}+2\epsilon},$
we have that
\begin{eqnarray}
&&K_{1}(\xi_{1},\mu_{1},\tau_{1},\xi,\mu,\tau)\leq\frac{|\xi|\langle\xi\rangle^{s_{1}}}
{\langle\sigma\rangle^{\frac{1}{2}-2\epsilon}
\prod\limits_{j=1}^{2}\langle\xi_{j}\rangle^{s_{1}}\langle\sigma_{j}\rangle^{\frac{1}{2}+\epsilon}}
\leq C\frac{|\xi|\langle\xi\rangle^{s_{1}}}
{\langle\sigma_{1}\rangle^{\frac{1}{2}-2\epsilon}\langle\sigma_{2}\rangle^{\frac{1}{2}+\epsilon}\langle\sigma\rangle^{\frac{1}{2}+\epsilon}
\prod\limits_{j=1}^{2}\langle\xi_{j}\rangle^{s_{1}}}
\nonumber\\
&&\leq C\frac{|\xi|^{s_{1}+\frac{1}{2}+2\epsilon}|\xi_{2}|^{-2s_{1}-\frac{\alpha}{2}+2\alpha\epsilon}}
{\langle\sigma_{2}\rangle^{\frac{1}{2}+\epsilon}\langle\sigma\rangle^{\frac{1}{2}+\epsilon}}
;\label{3.027}
\end{eqnarray}
when $s_{1}+\frac{1}{2}+2\epsilon\leq-\frac{1}{2},$ from  (\ref{3.027}), since  $s_{1}\geq\frac{1}{4}-\frac{3\alpha}{8}+4\alpha\epsilon$, we have that
\begin{eqnarray}
K_{1}(\xi_{1},\mu_{1},\tau_{1},\xi,\mu,\tau)\leq C\frac{|\xi|^{-\frac{1}{2}}|\xi_{2}|^{\frac{\alpha}{4}}}
{\langle\sigma\rangle^{\frac{1}{2}+\epsilon}\langle\sigma_{2}\rangle^{\frac{1}{2}+\epsilon}};\label{3.028}
\end{eqnarray}
when $s_{1}+\frac{1}{2}+2\epsilon>-\frac{1}{2},$ from  (\ref{3.027}), since  $s_{1}\geq\frac{1}{4}-\frac{3\alpha}{8}+4\alpha\epsilon$, we have that
\begin{eqnarray}
&&K_{1}(\xi_{1},\mu_{1},\tau_{1},\xi,\mu,\tau)\leq C\frac{|\xi|^{s_{1}+1+2\epsilon}|\xi|^{-\frac{1}{2}}|\xi_{1}|^{-2s_{1}-\frac{\alpha}{2}+2\alpha\epsilon}}
{\langle\sigma\rangle^{\frac{1}{2}+\epsilon}\langle\sigma_{2}\rangle^{\frac{1}{2}+\epsilon}}\nonumber\\&&\leq C\frac{|\xi|^{-\frac{1}{2}}|\xi_{2}|^{-s_{1}+1-\frac{\alpha}{2}+(2\alpha+2)\epsilon}}
{\langle\sigma\rangle^{\frac{1}{2}+\epsilon}\langle\sigma_{2}\rangle^{\frac{1}{2}+\epsilon}}\leq C\frac{|\xi|^{-\frac{1}{2}}|\xi_{2}|^{\frac{\alpha}{4}}}
{\langle\sigma\rangle^{\frac{1}{2}+\epsilon}\langle\sigma_{2}\rangle^{\frac{1}{2}+\epsilon}};\label{3.029}
\end{eqnarray}
Thus, combining (\ref{2.013}) with (\ref{3.028})-(\ref{3.029}),  we have that
\begin{eqnarray*}
|{\rm Int_{5}}|\leq C\|F\|_{L_{\tau\xi\mu}^{2}}\left(\prod\limits_{j=1}^{2}\|F_{j}\|_{L_{\tau\xi\mu}^{2}}\right).
\end{eqnarray*}
When (\ref{3.018}) is valid, this case can be proved similarly to (\ref{3.017}) of Region 5 with the aid of (\ref{2.016}).

\noindent When  $s_{1}\geq0,$ we have that
\begin{eqnarray}
&&K_{1}(\xi_{1},\mu_{1},\tau_{1},\xi,\mu,\tau)\leq C\frac{|\xi|}
{\langle\sigma\rangle^{\frac{1}{2}-2\epsilon}\prod\limits_{j=1}^{2}\langle\sigma_{j}\rangle^{\frac{1}{2}+\epsilon}}.\label{3.030}
\end{eqnarray}
We consider (\ref{3.016})-(\ref{3.018}),  respectively.

\noindent When (\ref{3.016}) is valid, from (\ref{3.030}), we have that
\begin{eqnarray}
&&K_{1}(\xi_{1},\mu_{1},\tau_{1},\xi,\mu,\tau)\leq C\frac{|\xi|^{\frac{1}{2}+2\epsilon}|\xi_{1}|^{-\frac{\alpha}{2}+2\alpha\epsilon}}
{\prod\limits_{j=1}^{2}\langle\sigma_{j}\rangle^{\frac{1}{2}+\epsilon}}\leq C\frac{\prod\limits_{j=1}^{2}|\xi_{j}|^{-\frac{\alpha}{8}+\frac{1}{4}}}
{\prod\limits_{j=1}^{2}\langle\sigma_{j}\rangle^{\frac{1}{2}+\epsilon}}.\label{3.031}
\end{eqnarray}
This case can be proved similarly to (\ref{3.025}).

\noindent When (\ref{3.017}) is valid, from (\ref{3.030}), since $\langle \sigma\rangle^{-\frac{1}{2}+2\epsilon}\langle \sigma_{1}\rangle^{-\frac{1}{2}-\epsilon}
 \leq \langle \sigma\rangle^{-\frac{1}{2}-\epsilon}\langle \sigma_{1}\rangle^{-\frac{1}{2}+2\epsilon}$ and $\alpha\geq4,$ we have that
\begin{eqnarray}
&&K_{1}(\xi_{1},\mu_{1},\tau_{1},\xi,\mu,\tau)\leq C\frac{|\xi|}{\langle\sigma_{1}\rangle^{\frac{1}{2}-2\epsilon}\langle\sigma\rangle^{\frac{1}{2}+\epsilon}\langle\sigma_{2}\rangle^{\frac{1}{2}+\epsilon}}\leq C\frac{|\xi|^{\frac{1}{2}+2\epsilon}|\xi_{1}|^{-\frac{\alpha}{2}+2\alpha\epsilon}}
{\langle\sigma\rangle^{\frac{1}{2}+\epsilon}\langle\sigma_{2}\rangle^{\frac{1}{2}+\epsilon}}\nonumber\\&&\leq C\frac{|\xi|^{-\frac{1}{2}}|\xi_{2}|^{\frac{\alpha}{4}}}
{\langle\sigma\rangle^{\frac{1}{2}+\epsilon}\langle\sigma_{2}\rangle^{\frac{1}{2}+\epsilon}}.\label{3.032}
\end{eqnarray}
This case can be proved similarly to (\ref{3.029}).

\noindent When (\ref{3.017}) is valid, from (\ref{3.030}), since $\langle \sigma\rangle^{-\frac{1}{2}+2\epsilon}\langle \sigma_{2}\rangle^{-\frac{1}{2}-\epsilon}
 \leq \langle \sigma\rangle^{-\frac{1}{2}-\epsilon}\langle \sigma_{2}\rangle^{-\frac{1}{2}+2\epsilon}$ and $\alpha\geq4,$ we have that
\begin{eqnarray}
&&K_{1}(\xi_{1},\mu_{1},\tau_{1},\xi,\mu,\tau)\leq C\frac{|\xi|}{\langle\sigma_{2}\rangle^{\frac{1}{2}-2\epsilon}\langle\sigma\rangle^{\frac{1}{2}+\epsilon}\langle\sigma_{1}\rangle^{\frac{1}{2}+\epsilon}}\leq C\frac{|\xi|^{\frac{1}{2}+2\epsilon}|\xi_{1}|^{-\frac{\alpha}{2}+2\alpha\epsilon}}
{\langle\sigma\rangle^{\frac{1}{2}+\epsilon}\langle\sigma_{1}\rangle^{\frac{1}{2}+\epsilon}}\nonumber\\&&\leq C\frac{|\xi|^{-\frac{1}{2}}|\xi_{2}|^{\frac{\alpha}{4}}}
{\langle\sigma\rangle^{\frac{1}{2}+\epsilon}\langle\sigma_{1}\rangle^{\frac{1}{2}+\epsilon}}.\label{3.033}
\end{eqnarray}
This case can be proved similarly to (\ref{3.029}).

\noindent(6). Region $\Omega_{6}.$

\noindent
In this region, we consider (\ref{3.016})-(\ref{3.018}),  respectively.

\noindent When (\ref{3.016}) is valid, since  $s_{1}\geq\frac{1}{4}-\frac{3\alpha}{8}+4\alpha\epsilon$ and $\alpha\geq4$, we have that
\begin{eqnarray}
&&K_{1}(\xi_{1},\mu_{1},\tau_{1},\xi,\mu,\tau)\leq C\frac{|\xi|\langle\xi\rangle^{s_{1}}}
{\langle\sigma\rangle^{\frac{1}{2}-2\epsilon}
\prod\limits_{j=1}^{2}\langle\xi_{j}\rangle^{s_{1}}\langle\sigma_{j}\rangle^{\frac{1}{2}+\epsilon}}\nonumber\\
&&\leq C\frac{|\xi_{2}|^{\frac{1-\alpha}{2}-s_{1}+2(\alpha+1)\epsilon}}
{\prod\limits_{j=1}^{2}\langle\sigma_{j}\rangle^{\frac{1}{2}+\epsilon}}
\leq
C\frac{|\xi_{1}|^{-\frac{1}{2}}|\xi_{2}|^{\frac{\alpha}{4}}}
{\prod\limits_{j=1}^{2}\langle\sigma_{j}\rangle^{\frac{1}{2}+\epsilon}}
.\label{3.034}
\end{eqnarray}
Thus, combining (\ref{2.011}) with (\ref{3.034}),  we have that
\begin{eqnarray*}
|{\rm Int_{6}}|\leq C\|F\|_{L_{\tau\xi\mu}^{2}}\left(\prod\limits_{j=1}^{2}\|F_{j}\|_{L_{\tau\xi\mu}^{2}}\right).
\end{eqnarray*}
When (\ref{3.016}) is valid, since  $s_{1}\geq\frac{1}{4}-\frac{3\alpha}{8}+4\alpha\epsilon$ and $\alpha\geq4$ and
$\langle \sigma\rangle^{-\frac{1}{2}+2\epsilon}\langle \sigma_{1}\rangle^{-\frac{1}{2}-\epsilon}
\leq \langle \sigma\rangle^{-\frac{1}{2}-\epsilon}\langle \sigma_{1}\rangle^{-\frac{1}{2}+2\epsilon},$
we have that
\begin{eqnarray}
&&K_{1}(\xi_{1},\mu_{1},\tau_{1},\xi,\mu,\tau)\leq C\frac{|\xi|\langle\xi\rangle^{s_{1}}}
{\langle\sigma\rangle^{\frac{1}{2}-2\epsilon}
\prod\limits_{j=1}^{2}\langle\xi_{j}\rangle^{s_{1}}\langle\sigma_{j}\rangle^{\frac{1}{2}+\epsilon}}\nonumber\\
&&\leq C\frac{|\xi_{2}|^{\frac{1-\alpha}{2}-s_{1}+2(\alpha+1)\epsilon}}
{\langle\sigma_{2}\rangle^{\frac{1}{2}+\epsilon}\langle\sigma\rangle^{\frac{1}{2}+\epsilon}}
\leq
C\frac{|\xi_{2}|^{-\frac{1}{2}}|\xi|^{\frac{\alpha}{4}}}
{\langle\sigma_{2}\rangle^{\frac{1}{2}+\epsilon}\langle\sigma\rangle^{\frac{1}{2}+\epsilon}}
.\label{3.035}
\end{eqnarray}
Thus, combining (\ref{2.015}) with (\ref{3.035}),  we have that
\begin{eqnarray*}
{\rm |Int_{6}|}\leq C\|F\|_{L_{\tau\xi\mu}^{2}}\left(\prod\limits_{j=1}^{2}\|F_{j}\|_{L_{\tau\xi\mu}^{2}}\right).
\end{eqnarray*}
When (\ref{3.018}) is valid, this case can be proved similarly to Region 6 of (\ref{3.017}) with the aid of (\ref{2.014}).

\noindent (7). Region $\Omega_{7}.$ This case can be proved similarly to Region $\Omega_{6}.$

This ends the proof of Lemma 3.1.

\noindent{\bf Remark 3.}
(\ref{3.011})  determines  $\alpha\geq4$ and
case (\ref{3.021})  of  Region 4 requires since  $s_{1}\geq\frac{1}{4}-\frac{3\alpha}{8}+4\alpha\epsilon$. When $s_{1}\geq\frac{1}{4}-\frac{3\alpha}{8}+4\alpha\epsilon$, we have that $-s_{1}<-\frac{1}{4}+\frac{3\alpha}{8}-4\alpha\epsilon.$

 \begin{Lemma}\label{Lemma3.2}
Let $-\frac{3\alpha}{8}+\frac{1}{2}+4\alpha\epsilon\leq s<0$ and $0<\epsilon<\frac{1}{100\alpha}$.
Then, we have that
\begin{eqnarray}
&&\hspace{-1cm}\|\partial_{x}\left[I_{N}(u_{1}u_{2})-I_{N}u_{1}I_{N}u_{2}\right]
\|_{X_{-\frac{1}{2}+2\epsilon}^{0,0}}\leq CN^{-\frac{3\alpha}{4}+1+(2\alpha+2)\epsilon}
\left(\prod_{j=1}^{2}\|I_{N}u_{j}\|_{X_{\frac{1}{2}+\epsilon}^{0,0}}\right).\label{3.036}
\end{eqnarray}
\end{Lemma}
\noindent{\bf Proof.}
To prove (\ref{3.036}),  by duality, it suffices to  prove that
\begin{eqnarray}
&&\left|\int_{\SR^{3}}\bar{h}\partial_{x}\left[I_{N}(u_{1}u_{2})-I_{N}u_{1}I_{N}u_{2}\right]dxdydt\right|\nonumber\\&&\leq
CN^{-\frac{3\alpha}{4}+1+(2\alpha+2)\epsilon}\|h\|_{X_{\frac{1}{2}-2\epsilon}^{0,0}}\left(\prod_{j=1}^{2}
\|I_{N}u_{j}\|_{X_{\frac{1}{2}+\epsilon}^{0,0}}\right)\label{3.037}
\end{eqnarray}
for $h\in X_{\frac{1}{2}-2\epsilon}^{0,0}.$
Let
\begin{eqnarray}
&&\xi=\xi_{1}+\xi_{2},\mu=\mu_{1}+\mu_{2},\tau=\tau_{1}+\tau_{2},\nonumber\\
&&F(\xi,\mu,\tau)=
\langle \sigma\rangle^{\frac{1}{2}-2\epsilon}M(\xi)\mathscr{F}h(\xi,\mu,\tau),\nonumber\\&&
F_{j}(\xi_{j},\mu_{j},\tau_{j})=M(\xi_{j})
\langle \sigma_{j}\rangle^{\frac{1}{2}+\epsilon}
\mathscr{F}u_{j}(\xi_{j},\mu,\tau_{j})(j=1,2).\label{3.038}
\end{eqnarray}
To obtain (\ref{3.037}), from (\ref{3.038}), it suffices to prove that
\begin{eqnarray}
&&\int_{D^{*}}\frac{|\xi|G(\xi_{1},\xi_{2})
F(\xi,\mu,\tau)\prod\limits_{j=1}^{2}F_{j}(\xi_{j},\mu_{j},\tau_{j})}{\langle\sigma_{j}\rangle^{\frac{1}{2}-2\epsilon}
\prod\limits_{j=1}^{2}\langle\sigma_{j}\rangle^{\frac{1}{2}+\epsilon}}
d\xi_{1}d\mu_{1}d\tau_{1}d\xi d\mu d\tau\nonumber\\&&\leq CN^{-\frac{3\alpha}{4}+1+(2\alpha+2)\epsilon}
\|F\|_{L_{\tau\xi\mu}^{2}}\left(\prod_{j=1}^{2}\|F_{j}\|_{L_{\tau\xi\mu}^{2}}\right),\label{3.039}
\end{eqnarray}
where
\begin{eqnarray*}
G(\xi_{1},\xi_{2})=\frac{M(\xi_{1})M(\xi_{2})-M(\xi)}{M(\xi_{1})M(\xi_{2})}
\end{eqnarray*}
and $D^{*}$ is defined as in Lemma 3.1.
Without loss of generality,  we assume that
 $F(\xi,\mu,\tau)\geq 0,F_j(\xi_{j},\mu_{j},\tau_{j})\geq 0(j=1,2)$.
By symmetry, we can assume that $|\xi_{1}|\geq |\xi_{2}|.$

 \noindent
We define
\begin{eqnarray*}
&&\hspace{-0.8cm}A_1=\left\{(\xi_1,\mu_{1},\tau_1,\xi,\mu,\tau)\in D^{*},
  |\xi_{2}|\leq |\xi_{1}|\leq \frac{N}{2}\right\},\\
&&\hspace{-0.8cm} A_2=\{ (\xi_1,\mu_{1},\tau_1,\xi,\mu,\tau)\in D^{*},
|\xi_1|>\frac{N}{2},|\xi_{1}|\geq |\xi_{2}|, |\xi_{2}|\leq 2\},\\
&&\hspace{-0.8cm} A_3=\{ (\xi_1,\mu_{1},\tau_1,\xi,\mu,\tau)\in D^{*},
|\xi_1|>\frac{N}{2},|\xi_{1}|\geq |\xi_{2}|, 2<|\xi_{2}|\leq N\},\\
&&\hspace{-0.8cm}A_4=\{(\xi_1,\mu_{1},\tau_1,\xi,\mu,\tau)\in D^{*},
|\xi_{1}|>\frac{N}{2},|\xi_{1}|\geq |\xi_{2}|,|\xi_{2}|>N\}.
\end{eqnarray*}
Here $D^{*}$ is defined as in Lemma 3.1.
Obviously, $D^{*}\subset\bigcup\limits_{j=1}^{4}A_{j}.$
We define
\begin{equation}
    K_{2}(\xi_{1},\mu_{1},\tau_{1},\xi,\mu,\tau):=\frac{|\xi|G(\xi_{1},\xi_{2})}
    {\langle\sigma_{j}\rangle^{\frac{1}{2}-2\epsilon}
\prod\limits_{j=1}^{2}\langle\sigma_{j}\rangle^{\frac{1}{2}+\epsilon}}\label{3.040}
\end{equation}
and
\begin{eqnarray*}
J_{k}:=\int_{A_{k}} K_{2}(\xi_{1},\mu_{1},\tau_{1},\xi,\mu,\tau)F(\xi,\mu,\tau)
\prod_{j=1}^{2}F_{j}(\xi_{j},\mu_{j},\tau_{j})
d\xi_{1}d\mu_{1}d\tau_{1}d\xi d\mu d\tau
\end{eqnarray*}
with $1\leq k\leq 4, k\in N.$

\noindent Since (\ref{2.024}) is valid,
one of (\ref{3.016})-(\ref{3.018}) must occur,

\noindent (1) Region $A_{1}.$ In this case, since $G(\xi_{1},\xi_{2})=0$, thus we have that $J_{1}=0$.

\noindent (2) Region $A_{2}$. From page 902 of \cite{ILM-CPAA}, we have that
\begin{eqnarray}
G(\xi_{1},\xi_{2})\leq C\frac{|\xi_{2}|}{|\xi_{1}|}\label{3.041}.
\end{eqnarray}
Inserting (\ref{3.041}) into (\ref{3.040}) yields
\begin{eqnarray}
    K_{2}(\xi_{1},\mu_{1},\tau_{1},\xi,\mu,\tau)\leq C\frac{|\xi|G(\xi_{1},\xi_{2})}
    {\langle\sigma\rangle^{\frac{1}{2}-2\epsilon}
\prod\limits_{j=1}^{2}\langle\sigma_{j}\rangle^{\frac{1}{2}+\epsilon}}\leq
\frac{C|\xi_{2}|}{\langle\sigma\rangle^{\frac{1}{2}-2\epsilon}
\prod\limits_{j=1}^{2}\langle\sigma_{j}\rangle^{\frac{1}{2}+\epsilon}}\label{3.042}.
\end{eqnarray}
When (\ref{3.016}) is valid,  since $0<\epsilon<\frac{1}{100\alpha}$,  then  $\frac{1}{2}-2\epsilon>0,$ thus, we have that
\begin{eqnarray}
    K_{2}(\xi_{1},\mu_{1},\tau_{1},\xi,\mu,\tau)\leq
\frac{C|\xi_{2}|^{\frac{1}{2}+2\epsilon}|\xi_{1}|^{-\frac{\alpha}{2}+2\alpha\epsilon}}{
\prod\limits_{j=1}^{2}\langle\sigma_{j}\rangle^{\frac{1}{2}+\epsilon}}\leq
\frac{CN^{-\frac{3\alpha}{4}+2\alpha\epsilon}|\xi_{2}|^{-\frac{1}{2}}|\xi_{1}|^{\frac{\alpha}{4}}}{
\prod\limits_{j=1}^{2}\langle\sigma_{j}\rangle^{\frac{1}{2}+\epsilon}}\label{3.043}.
\end{eqnarray}
Using (\ref{2.011}) and (\ref{3.043}),  we have that
\begin{eqnarray}
&&J_{2}\leq C N^{-\frac{3\alpha}{4}+2\alpha\epsilon}\int_{A_{2}}\frac{|\xi_{2}|^{-\frac{1}{2}}|\xi_{1}|^{\frac{\alpha}{4}}}{
\prod\limits_{j=1}^{2}\langle\sigma_{j}\rangle^{\frac{1}{2}+\epsilon}}F(\xi,\mu,\tau)
\prod_{j=1}^{2}F_{j}(\xi_{j},\mu_{j},\tau_{j})d\xi_{1}d\mu_{1}d\tau_{1}d\xi d\mu d\tau\nonumber\\
&&\leq CN^{-\frac{3\alpha}{4}+2\alpha\epsilon}\|F\|_{L_{\tau\xi\mu}^{2}}\left(\prod_{j=1}^{2}\|F_{j}\|_{L_{\tau\xi\mu}^{2}}\right)\label{3.044}.
\end{eqnarray}
When (\ref{3.017}) is valid, since $\langle \sigma\rangle^{-\frac{1}{2}+2\epsilon}\langle \sigma_{1}\rangle^{-\frac{1}{2}-\epsilon}
\leq \langle \sigma\rangle^{-\frac{1}{2}-\epsilon}\langle \sigma_{1}\rangle^{-\frac{1}{2}+2\epsilon}$ and  $0<\epsilon<\frac{1}{100\alpha}$, we have that
\begin{eqnarray}
    K_{2}(\xi_{1},\mu_{1},\tau_{1},\xi,\mu,\tau)\leq
\frac{C|\xi_{2}|^{\frac{1}{2}+2\epsilon}|\xi_{1}|^{-\frac{\alpha}{2}+2\alpha\epsilon}}{
\langle\sigma_{2}\rangle^{\frac{1}{2}+\epsilon}\langle\sigma\rangle^{\frac{1}{2}+\epsilon}}\leq
\frac{CN^{-\frac{3\alpha}{4}+2\alpha\epsilon}|\xi_{2}|^{-\frac{1}{2}}|\xi|^{\frac{\alpha}{4}}}{
\langle\sigma_{2}\rangle^{\frac{1}{2}+\epsilon}\langle\sigma\rangle^{\frac{1}{2}+\epsilon}}\label{3.045}.
\end{eqnarray}
By using (\ref{2.015}), from (\ref{3.045}),  we have that
\begin{eqnarray}
&&J_{2}\leq C N^{-\frac{3\alpha}{4}+2\alpha\epsilon}\int_{A_{2}}\frac{|\xi_{2}|^{-\frac{1}{2}}|\xi|^{\frac{\alpha}{4}}}{
\langle\sigma_{2}\rangle^{\frac{1}{2}+\epsilon}\langle\sigma\rangle^{\frac{1}{2}+\epsilon}}F(\xi,\mu,\tau)
\prod_{j=1}^{2}F_{j}(\xi_{j},\mu_{j},\tau_{j})d\xi_{1}d\mu_{1}d\tau_{1}d\xi d\mu d\tau\nonumber\\
&&\leq CN^{-\frac{3\alpha}{4}+2\alpha\epsilon}\|F\|_{L_{\tau\xi\mu}^{2}}\left(\prod_{j=1}^{2}\|F_{j}\|_{L_{\tau\xi\mu}^{2}}\right)\label{3.046}.
\end{eqnarray}
When (\ref{3.018}) is valid,  since $\langle \sigma\rangle^{-\frac{1}{2}+2\epsilon}\langle \sigma_{2}\rangle^{-\frac{1}{2}-\epsilon}
\leq \langle \sigma\rangle^{-\frac{1}{2}-\epsilon}\langle \sigma_{2}\rangle^{-\frac{1}{2}+2\epsilon}$ and  $0<\epsilon<\frac{1}{100\alpha}$, we have that
\begin{eqnarray}
    K_{2}(\xi_{1},\mu_{1},\tau_{1},\xi,\mu,\tau)\leq
\frac{C|\xi_{2}|^{\frac{1}{2}+2\epsilon}|\xi_{1}|^{-\frac{\alpha}{2}+2\alpha\epsilon}}{
\langle\sigma_{1}\rangle^{\frac{1}{2}+\epsilon}\langle\sigma\rangle^{\frac{1}{2}+\epsilon}}\leq
\frac{CN^{-\frac{3\alpha}{4}+\frac{1}{2}+2\alpha\epsilon}|\xi_{1}|^{-\frac{1}{2}}|\xi|^{\frac{\alpha}{4}}}{
\langle\sigma_{1}\rangle^{\frac{1}{2}+\epsilon}\langle\sigma\rangle^{\frac{1}{2}+\epsilon}}\label{3.047}.
\end{eqnarray}
By using (\ref{2.012}), from (\ref{3.047}),  we have that
\begin{eqnarray}
&&J_{2}\leq C N^{-\frac{3\alpha}{4}+\frac{1}{2}+2\alpha\epsilon}\int_{A_{2}}\frac{|\xi_{1}|^{-\frac{1}{2}}|\xi|^{\frac{\alpha}{4}}}{
\langle\sigma_{2}\rangle^{\frac{1}{2}+\epsilon}\langle\sigma\rangle^{\frac{1}{2}+\epsilon}}F(\xi,\mu,\tau)
\prod_{j=1}^{2}F_{j}(\xi_{j},\mu_{j},\tau_{j})d\xi_{1}d\mu_{1}d\tau_{1}d\xi d\mu d\tau\nonumber\\
&&\leq CN^{-\frac{3\alpha}{4}+\frac{1}{2}+2\alpha\epsilon}\|F\|_{L_{\tau\xi\mu}^{2}}\left(\prod_{j=1}^{2}\|F_{j}\|_{L_{\tau\xi\mu}^{2}}\right)\label{3.048}.
\end{eqnarray}
(3) Region $A_{3}$.
From page 902 of \cite{ILM-CPAA}, we have that (\ref{3.041}) is valid.
Combining (\ref{3.041}) with (\ref{3.040}), we have that
\begin{equation}
    K_{2}(\xi_{1},\mu_{1},\tau_{1},\xi,\mu,\tau)\leq C\frac{{\rm min}
    \left\{|\xi|,|\xi_{1}|,|\xi_{2}|\right\}}{\langle\sigma\rangle^{\frac{1}{2}-2\epsilon}
\prod\limits_{j=1}^{2}\langle\sigma_{j}\rangle^{\frac{1}{2}+\epsilon}}.\label{3.049}
\end{equation}
In this case, we consider $|\xi_{1}|>6|\xi_{2}|$, $|\xi_{1}|\leq 6|\xi_{2}|,$ respectively.

\noindent When $|\xi_{1}|> 6|\xi_{2}|$, we consider (\ref{3.016})-(\ref{3.018}), respectively.

\noindent When (\ref{3.016})  is valid, from (\ref{3.049}), since $0<\epsilon<\frac{1}{100\alpha}$, we have that
\begin{eqnarray}
   && K_{2}(\xi_{1},\mu_1,\tau_{1},\xi,\mu,\tau)\leq C\frac{|\xi_{2}|^{\frac{1}{2}+2\epsilon}|\xi_{1}|^{-\frac{\alpha}{2}+2\alpha\epsilon}}{
\prod\limits_{j=1}^{2}\langle\sigma_{j}\rangle^{\frac{1}{2}+\epsilon}}\leq C\frac{|\xi_{2}|^{-\frac{1}{2}}|\xi_{2}|^{1+2\epsilon}|\xi_{1}|^{-\frac{\alpha}{2}+2\alpha\epsilon}}{\prod\limits_{j=1}^{2}\langle\sigma_{j}\rangle^{\frac{1}{2}+\epsilon}}
\nonumber\\&&\leq C
\frac{|\xi_{2}|^{-\frac{1}{2}}|\xi_{1}|^{-\frac{\alpha}{2}+1+(2\alpha+2)\epsilon}}{\prod\limits_{j=1}^{2}\langle\sigma_{j}\rangle^{\frac{1}{2}+\epsilon}}\leq C N^{-\frac{3\alpha}{4}+1+(2\alpha+2)\epsilon}\frac{|\xi_{2}|^{-\frac{1}{2}}|\xi_{1}|^{\frac{\alpha}{4}}}
{\prod\limits_{j=1}^{2}\langle\sigma_{j}\rangle^{\frac{1}{2}+\epsilon}}.\label{3.050}
\end{eqnarray}
Using   (\ref{2.011}),  from  (\ref{3.050}),  we
have that
\begin{eqnarray*}
&&J_{3}\leq CN^{-\frac{3\alpha}{4}+1+(2\alpha+2)\epsilon}\int_{A_{3}}\frac{|\xi_{2}|^{-\frac{1}{2}}|\xi_{1}|^{\frac{\alpha}{4}}F\prod\limits_{j=1}^{2}F_{j}}
{\prod\limits_{j=1}^{2}\langle\sigma_{j}\rangle^{\frac{1}{2}+\epsilon}}d\xi_{1}d\mu_{1}d\tau_{1}d\xi d\mu d\tau\nonumber\\
&&\leq CN^{-\frac{3\alpha}{4}+1+(2\alpha+2)\epsilon}
\|F\|_{L_{\tau\xi\mu}^{2}}\left(\prod_{j=1}^{2}\|F_{j}\|_{L_{\tau\xi\mu}^{2}}\right).
\end{eqnarray*}
When (\ref{3.017}) is valid, since $\langle \sigma\rangle ^{-\frac{1}{2}+2\epsilon}\langle \sigma_{1}\rangle^{-\frac{1}{2}-\epsilon}\leq \langle \sigma_{1}\rangle ^{-\frac{1}{2}+2\epsilon}\langle \sigma\rangle^{-\frac{1}{2}-\epsilon}$ and $0<\epsilon<\frac{1}{100\alpha}$, we have that
\begin{eqnarray}
   && K_{2}(\xi_{1},\mu_1,\tau_{1},\xi,\mu,\tau)\leq C\frac{|\xi_{2}|}{\langle\sigma_{1}\rangle^{\frac{1}{2}-2\epsilon}
   \langle\sigma_{2}\rangle^{\frac{1}{2}+\epsilon}\langle\sigma\rangle^{\frac{1}{2}+\epsilon}}\nonumber\\&&\leq C\frac{|\xi_{2}|^{\frac{1}{2}+2\epsilon}|\xi_{1}|^{-\frac{\alpha}{2}+2\alpha\epsilon}}{\langle\sigma_{2}\rangle^{\frac{1}{2}+\epsilon}
   \langle\sigma\rangle^{\frac{1}{2}+\epsilon}}
   \leq C
N^{-\frac{3\alpha}{4}+1+(2\alpha+2)\epsilon}
\frac{|\xi_{2}|^{-\frac{1}{2}}|\xi|^{\frac{\alpha}{4}}}{\langle\sigma_{2}\rangle^{\frac{1}{2}+\epsilon}
\langle\sigma\rangle^{\frac{1}{2}+\epsilon}}.\label{3.051}
\end{eqnarray}
Combining (\ref{2.013}) with (\ref{3.051}), we have that
\begin{eqnarray*}
&&J_{3}\leq CN^{-\frac{3\alpha}{4}+1+(2\alpha+2)\epsilon}\int_{A_{3}}\frac{|\xi_{2}|^{-\frac{1}{2}}
|\xi|^{\frac{\alpha}{4}}F\prod\limits_{j=1}^{2}F_{j}}{\langle\sigma_{2}\rangle^{\frac{1}{2}+\epsilon}
\langle\sigma\rangle^{\frac{1}{2}+\epsilon}}d\xi_{1}d\mu_{1}d\tau_{1}d\xi d\mu d\tau\nonumber\\
&&\leq  CN^{-\frac{3\alpha}{4}+1+(2\alpha+2)\epsilon}
\|F\|_{L_{\tau\xi\mu}^{2}}\left(\prod_{j=1}^{2}\|F_{j}\|_{L_{\tau\xi\mu}^{2}}\right).
\end{eqnarray*}
When (\ref{3.018}) is valid,  since $\langle \sigma\rangle ^{-\frac{1}{2}+2\epsilon}\langle \sigma_{2}\rangle^{-\frac{1}{2}-\epsilon}\leq \langle \sigma_{2}\rangle ^{-\frac{1}{2}+2\epsilon}\langle \sigma\rangle^{-\frac{1}{2}-\epsilon}$ and  $0<\epsilon<\frac{1}{100\alpha}$,
we have that
\begin{eqnarray}
   && K_{2}(\xi_{1},\mu_1,\tau_{1},\xi,\mu,\tau)\leq C\frac{|\xi_{2}|}{\langle\sigma_{1}\rangle^{\frac{1}{2}-2\epsilon}
   \langle\sigma_{2}\rangle^{\frac{1}{2}+\epsilon}\langle\sigma\rangle^{\frac{1}{2}+\epsilon}}\nonumber\\&&\leq C\frac{|\xi_{2}|^{\frac{1}{2}+2\epsilon}|\xi_{1}|^{-\frac{\alpha}{2}+2\alpha\epsilon}}{\langle\sigma_{1}\rangle^{\frac{1}{2}+\epsilon}
   \langle\sigma\rangle^{\frac{1}{2}+\epsilon}}
   \leq C
N^{-\frac{3\alpha}{4}+1+(2\alpha+2)\epsilon}
\frac{|\xi_{1}|^{-\frac{1}{4}+\frac{\alpha}{8}}|\xi|^{-\frac{1}{4}+\frac{\alpha}{8}}}{\langle\sigma_{1}\rangle^{\frac{1}{2}+\epsilon}
\langle\sigma\rangle^{\frac{1}{2}+\epsilon}}.\label{3.052}
\end{eqnarray}
Combining (\ref{2.09})  with (\ref{3.052}),  we have that
\begin{eqnarray*}
&&J_{3}\leq CN^{-\frac{3\alpha}{4}+1+(2\alpha+2)\epsilon}
\|F\|_{L_{\tau\xi\mu}^{2}}\left(\prod_{j=1}^{2}\|F_{j}\|_{L_{\tau\xi\mu}^{2}}\right).
\end{eqnarray*}
When $\frac{|\xi_{1}|}{6}\leq |\xi_{2}|\leq |\xi_{1}|,$ we consider (\ref{3.016})-(\ref{3.018}), respectively.

\noindent When (\ref{3.016})  is valid, from (\ref{3.049}),  we have that
\begin{equation}
    K_{2}(\xi_{1},\mu_1,\tau_{1},\xi,\mu,\tau)\leq C\frac{|\xi|^{\frac{1}{2}+2\epsilon}|\xi_{1}|^{-\frac{\alpha}{2}+2\alpha\epsilon}}{
\prod\limits_{j=1}^{2}\langle\sigma_{j}\rangle^{\frac{1}{2}+\epsilon}}\leq CN^{-\frac{3\alpha}{4}+1+(2\alpha+2)\epsilon}\frac{\prod\limits_{j=1}^{2}|\xi_{j}|^{-\frac{1}{4}+\frac{\alpha}{8}}}
{\prod\limits_{j=1}^{2}\langle\sigma_{j}\rangle^{\frac{1}{2}+\epsilon}}.\label{3.053}
\end{equation}
Combining (\ref{2.08}) with (\ref{3.053}),   we
have that
\begin{eqnarray*}
&&J_{3}\leq CN^{-\frac{3\alpha}{4}+1+(2\alpha+2)\epsilon}\int_{A_{3}}\frac{\prod\limits_{j=1}^{2}|\xi_{j}|^{-\frac{1}{4}+\frac{\alpha}{8}}F\prod\limits_{j=1}^{2}F_{j}}
{\prod\limits_{j=1}^{2}\langle\sigma_{j}\rangle^{\frac{1}{2}+\epsilon}}d\xi_{1}d\mu_{1}d\tau_{1}d\xi d\mu d\tau\nonumber\\
&&\leq CN^{-\frac{3\alpha}{4}+1+(2\alpha+2)\epsilon}
\|F\|_{L_{\tau\xi\mu}^{2}}\left(\prod_{j=1}^{2}\|F_{j}\|_{L_{\tau\xi\mu}^{2}}\right).
\end{eqnarray*}
When (\ref{3.017}) is valid,
since $\langle \sigma\rangle ^{-\frac{1}{2}+2\epsilon}\langle \sigma_{1}\rangle^{-\frac{1}{2}-\epsilon}\leq \langle \sigma_{1}\rangle ^{-\frac{1}{2}+2\epsilon}\langle \sigma\rangle^{-\frac{1}{2}-\epsilon}$, we have that
\begin{eqnarray}
  &&K_{2}(\xi_{1},\mu_1,\tau_{1},\xi,\mu,\tau)\leq C\frac{|\xi|}{\langle\sigma_{1}\rangle^{\frac{1}{2}-2\epsilon}
   \langle\sigma_{2}\rangle^{\frac{1}{2}+\epsilon}\langle\sigma\rangle^{\frac{1}{2}+\epsilon}}\nonumber\\&&\leq C \frac{|\xi|^{\frac{1}{2}+2\epsilon}|\xi_{1}|^{-\frac{\alpha}{2}+2\alpha\epsilon}}{\langle\sigma_{2}\rangle^{\frac{1}{2}+\epsilon}
   \langle\sigma\rangle^{\frac{1}{2}+\epsilon}}\leq CN^{-\frac{3\alpha}{4}+1+(2\alpha+2)\epsilon}\frac{|\xi|^{-\frac{1}{2}}|\xi_{2}|^{\frac{\alpha}{4}}}{\langle\sigma_{2}\rangle^{\frac{1}{2}+\epsilon}
   \langle\sigma\rangle^{\frac{1}{2}+\epsilon}}.
  \label{3.054}
\end{eqnarray}
Combining  (\ref{2.013}) with (\ref{3.054}),   we
have that
\begin{eqnarray*}
&&J_{3}\leq CN^{-\frac{3\alpha}{4}+1+(2\alpha+2)\epsilon}\int_{A_{3}}\frac{|\xi|^{-\frac{1}{2}}|\xi_{2}|^{\frac{\alpha}{4}}F\prod\limits_{j=1}^{2}F_{j}}
{\langle\sigma_{2}\rangle^{\frac{1}{2}+\epsilon}
   \langle\sigma\rangle^{\frac{1}{2}+\epsilon}}
   d\xi_{1}d\mu_{1}d\tau_{1}d\xi d\mu d\tau\nonumber\\
&&\leq CN^{-\frac{3\alpha}{4}+1+(2\alpha+2)\epsilon}
\|F\|_{L_{\tau\xi\mu}^{2}}\left(\prod_{j=1}^{2}\|F_{j}\|_{L_{\tau\xi\mu}^{2}}\right).
\end{eqnarray*}
When (\ref{3.018}) is valid,
this case can be proved similarly to case (\ref{3.017}).

\noindent (4) Region $A_{4}$.
From lines 19-20 of page 903 in \cite{ILM-CPAA}, we have that
\begin{eqnarray}
G(\xi_{1},\xi_{2})\leq C\prod_{j=1}^{2}\left(\frac{|\xi_{j}|}{N}\right)^{-s}.\label{3.055}
\end{eqnarray}
Inserting (\ref{3.055})  into  (\ref{3.040}) yields
\begin{eqnarray}
K_{2}(\xi_{1},\mu_1,\tau_{1},\xi,\mu,\tau)\leq CN^{2s}\frac{|\xi|\prod\limits_{j=1}^{2}|\xi_{j}|^{-s}}{
\langle\sigma\rangle^{\frac{1}{2}-2\epsilon}\prod\limits_{j=1}^{2}\langle\sigma_{j}\rangle^{\frac{1}{2}+\epsilon}}.\label{3.056}
\end{eqnarray}
We consider $|\xi|\leq\frac{N}{4}$, $|\xi|>\frac{N}{4}$,
respectively.

\noindent When $|\xi|\leq\frac{N}{4}$ is valid, we have that $|\xi_{1}|\sim |\xi_{2}|$,
we consider (\ref{3.016})-(\ref{3.018}),  respectively.

\noindent When (\ref{3.016}) is valid, from (\ref{3.056}), since $\frac{1}{2}-\frac{3\alpha}{8}+4\alpha\epsilon\leq s\leq0$ and $0<\epsilon<\frac{1}{100\alpha},$ we have that
\begin{eqnarray}
&&K_{2}(\xi_{1},\mu_1,\tau_{1},\xi,\mu,\tau)\leq CN^{2s}\frac{|\xi|^{\frac{1}{2}+2\epsilon}|\xi_{1}|^{-2s-\frac{\alpha}{2}+2\alpha\epsilon}}
{\prod\limits_{j=1}^{2}\langle\sigma_{j}\rangle^{\frac{1}{2}+\epsilon}}\leq CN^{2s}\frac{|\xi_{1}|^{-2s-\frac{\alpha}{2}+\frac{1}{2}+(2\alpha+2)\epsilon}}
{\prod\limits_{j=1}^{2}\langle\sigma_{j}\rangle^{\frac{1}{2}+\epsilon}}\nonumber\\
&&\leq CN^{2s}\frac{|\xi_{1}|^{-2s-\frac{3\alpha}{4}+1+(2\alpha+2)\epsilon}\prod\limits_{j=1}^{2}|\xi_{j}|^{-\frac{1}{4}+\frac{\alpha}{8}}}
{\prod\limits_{j=1}^{2}\langle\sigma_{j}\rangle^{\frac{1}{2}+\epsilon}}
\leq CN^{-\frac{3\alpha}{4}+1+(2\alpha+2)\epsilon}\frac{\prod\limits_{j=1}^{2}|\xi_{j}|^{-\frac{1}{4}+\frac{\alpha}{8}}}
{\prod\limits_{j=1}^{2}\langle\sigma_{j}\rangle^{\frac{1}{2}+\epsilon}}.\label{3.057}
\end{eqnarray}
Combining (\ref{2.08}) with (\ref{3.057}),  we have that
\begin{eqnarray*}
&&J_{4}\leq CN^{-\frac{3\alpha}{4}+1+(2\alpha+2)\epsilon}\int_{A_{4}}\frac{\prod\limits_{j=1}^{2}|\xi_{j}|^{-\frac{1}{4}+\frac{\alpha}{8}}F\prod\limits_{j=1}^{2}F_{j}}
{\langle\sigma_{2}\rangle^{\frac{1}{2}+\epsilon}
   \langle\sigma\rangle^{\frac{1}{2}+\epsilon}}
   d\xi_{1}d\mu_{1}d\tau_{1}d\xi d\mu d\tau\nonumber\\
&&\leq CN^{-\frac{3\alpha}{4}+1+(2\alpha+2)\epsilon}
\|F\|_{L_{\tau\xi\mu}^{2}}\left(\prod_{j=1}^{2}\|F_{j}\|_{L_{\tau\xi\mu}^{2}}\right).
\end{eqnarray*}
When (\ref{3.017}) is valid,
since $\langle \sigma\rangle ^{-\frac{1}{2}+2\epsilon}\langle \sigma_{1}\rangle^{-\frac{1}{2}-\epsilon}\leq \langle \sigma_{1}\rangle ^{-\frac{1}{2}+2\epsilon}\langle \sigma\rangle^{-\frac{1}{2}-\epsilon}$ and   $-\frac{3\alpha}{8}+\frac{1}{2}+4\alpha\epsilon\leq s\leq0$ and $|\xi_{1}|\sim |\xi_{2}|,$ we have that
\begin{eqnarray}
   && K_{2}(\xi_{1},\mu_1,\tau_{1},\xi,\mu,\tau)\leq CN^{2s}\frac{|\xi|^{\frac{1}{2}+2\epsilon}|\xi_{1}|^{-2s-\frac{\alpha}{2}+2\alpha\epsilon}}
{\langle\sigma_{2}\rangle^{\frac{1}{2}+\epsilon}\langle\sigma\rangle^{\frac{1}{2}+\epsilon}}\nonumber\\
&&\leq CN^{2s}\frac{|\xi|^{-\frac{1}{2}}|\xi_{1}|^{-2s-\frac{\alpha}{2}+1+(2\alpha+2)\epsilon}}
{\langle\sigma_{2}\rangle^{\frac{1}{2}+\epsilon}\langle\sigma\rangle^{\frac{1}{2}+\epsilon}}\nonumber\\
&&\leq CN^{2s}\frac{|\xi_{1}|^{-2s-\frac{3\alpha}{4}+1+(2\alpha+2)\epsilon}|\xi|^{-\frac{1}{2}}|\xi_{2}|^{\frac{\alpha}{4}}}
{\langle\sigma_{2}\rangle^{\frac{1}{2}+\epsilon}\langle\sigma\rangle^{\frac{1}{2}+\epsilon}}
\leq CN^{-\frac{3\alpha}{4}+1+(2\alpha+2)\epsilon}\frac{|\xi|^{-\frac{1}{2}}|\xi_{2}|^{\frac{\alpha}{4}}}
{\langle\sigma_{2}\rangle^{\frac{1}{2}+\epsilon}\langle\sigma\rangle^{\frac{1}{2}+\epsilon}}.
  \label{3.058}
\end{eqnarray}
Combining (\ref{2.013}) with (\ref{3.058}),  we have that
\begin{eqnarray*}
&&J_{4}\leq CN^{-\frac{3\alpha}{4}+1+(2\alpha+2)\epsilon}\int_{A_{4}}\frac{|\xi_{2}|^{-\frac{1}{2}}|\xi|^{\frac{\alpha}{4}}F\prod\limits_{j=1}^{2}F_{j}}
{\langle\sigma_{2}\rangle^{\frac{1}{2}+\epsilon}
   \langle\sigma\rangle^{\frac{1}{2}+\epsilon}}
   d\xi_{1}d\mu_{1}d\tau_{1}d\xi d\mu d\tau\nonumber\\
&&\leq CN^{-\frac{3\alpha}{4}+1+(2\alpha+2)\epsilon}
\|F\|_{L_{\tau\xi\mu}^{2}}\left(\prod_{j=1}^{2}\|F_{j}\|_{L_{\tau\xi\mu}^{2}}\right).
\end{eqnarray*}
 When (\ref{3.018}) is valid,
this case can be proved similarly to case  (\ref{3.017}).

\noindent When $|\xi|\geq \frac{N}{4},$ we consider $|\xi|<\frac{|\xi_{2}|}{6},|\xi|\geq \frac{|\xi_{2}|}{6},$ respectively.

\noindent When  $|\xi|<\frac{|\xi_{2}|}{6},$ we have that  $|\xi_{1}|\sim |\xi_{2}|,$ this case can be proved similarly to case $|\xi|<\frac{N}{4}.$

\noindent When $|\xi|\geq \frac{|\xi_{2}|}{6},$ we consider $\frac{|\xi_{2}|}{6}\leq |\xi|\leq 6|\xi_{2}|,|\xi|>6|\xi_{2}|,$ respectively.

\noindent When $\frac{|\xi_{2}|}{6}\leq |\xi|\leq 6|\xi_{2}|,$ we have that  $|\xi_{1}|\sim |\xi_{2}|\sim |\xi|,$ we consider (\ref{3.016})-(\ref{3.018}),  respectively.

\noindent When (\ref{3.016}) is valid, from (\ref{3.056}), since $-\frac{3\alpha}{8}+\frac{1}{2}+4\alpha\epsilon\leq s\leq0,$ we have that
\begin{eqnarray}
&&K_{2}(\xi_{1},\mu_1,\tau_{1},\xi,\mu,\tau)\leq CN^{2s}\frac{|\xi|^{\frac{1}{2}+2\epsilon}|\xi_{1}|^{-2s-\frac{\alpha}{2}+2\alpha\epsilon}}
{\prod\limits_{j=1}^{2}\langle\sigma_{j}\rangle^{\frac{1}{2}+\epsilon}}\nonumber\\
&&\leq CN^{2s}\frac{|\xi_{1}|^{-2s-\frac{3\alpha}{4}+\frac{1}{2}+(2\alpha+2)\epsilon}\prod\limits_{j=1}^{2}|\xi_{j}|^{-\frac{1}{4}+\frac{\alpha}{8}}}
{\prod\limits_{j=1}^{2}\langle\sigma_{j}\rangle^{\frac{1}{2}+\epsilon}}\leq CN^{-\frac{3\alpha}{4}+\frac{1}{2}+(2\alpha+2)\epsilon}\frac{\prod\limits_{j=1}^{2}|\xi_{j}|^{-\frac{1}{4}+\frac{\alpha}{8}}}
{\prod\limits_{j=1}^{2}\langle\sigma_{j}\rangle^{\frac{1}{2}+\epsilon}}.\label{3.059}
\end{eqnarray}
Combining (\ref{2.08}) with (\ref{3.059}),  we have that
\begin{eqnarray*}
&&J_{4}\leq CN^{-\frac{3\alpha}{4}+1+(2\alpha+2)\epsilon}\int_{A_{4}}\frac{\prod\limits_{j=1}^{2}|\xi_{j}|^{-\frac{1}{4}+\frac{\alpha}{8}}F\prod\limits_{j=1}^{2}F_{j}}
{\prod\limits_{j=1}^{2}\langle\sigma_{j}\rangle^{\frac{1}{2}+\epsilon}}
   d\xi_{1}d\mu_{1}d\tau_{1}d\xi d\mu d\tau\nonumber\\
&&\leq CN^{-\frac{3\alpha}{4}+1+(2\alpha+2)\epsilon}
\|F\|_{L_{\tau\xi\mu}^{2}}\left(\prod_{j=1}^{2}\|F_{j}\|_{L_{\tau\xi\mu}^{2}}\right).
\end{eqnarray*}
When (\ref{3.017}) is valid,
since $\langle \sigma\rangle ^{-\frac{1}{2}+2\epsilon}\langle \sigma_{1}\rangle^{-\frac{1}{2}-\epsilon}\leq \langle \sigma_{1}\rangle ^{-\frac{1}{2}+2\epsilon}\langle \sigma\rangle^{-\frac{1}{2}-\epsilon}$ and    $-\frac{3\alpha}{8}+\frac{1}{2}+4\alpha\epsilon\leq s\leq0,$ we have that
\begin{eqnarray}
&& K_{2}(\xi_{1},\mu_1,\tau_{1},\xi,\mu,\tau)\leq CN^{2s}\frac{|\xi|^{\frac{1}{2}+2\epsilon}|\xi_{1}|^{-2s-\frac{\alpha}{2}+2\alpha\epsilon}}
{\langle\sigma_{2}\rangle^{\frac{1}{2}+\epsilon}\langle\sigma\rangle^{\frac{1}{2}+\epsilon}}\nonumber\\
&&\leq CN^{2s}\frac{|\xi|^{-\frac{1}{2}}|\xi_{1}|^{-2s-\frac{\alpha}{2}+1+(2\alpha+2)\epsilon}}
{\langle\sigma_{2}\rangle^{\frac{1}{2}+\epsilon}\langle\sigma\rangle^{\frac{1}{2}+\epsilon}}\nonumber\\
&&\leq CN^{2s}\frac{|\xi_{1}|^{-2s-\frac{3\alpha}{4}+1+(2\alpha+2)\epsilon}|\xi|^{-\frac{1}{2}}|\xi_{2}|^{\frac{\alpha}{4}}}
{\langle\sigma_{2}\rangle^{\frac{1}{2}+\epsilon}\langle\sigma\rangle^{\frac{1}{2}+\epsilon}}
\leq CN^{-\frac{3\alpha}{4}+1+(2\alpha+2)\epsilon}\frac{|\xi|^{-\frac{1}{2}}|\xi_{2}|^{\frac{\alpha}{4}}}
{\langle\sigma_{2}\rangle^{\frac{1}{2}+\epsilon}\langle\sigma\rangle^{\frac{1}{2}+\epsilon}}.
  \label{3.060}
\end{eqnarray}
Combining (\ref{2.08}) with (\ref{3.060}),  we have that
\begin{eqnarray*}
&&J_{4}\leq CN^{-\frac{3\alpha}{4}+1+(2\alpha+2)\epsilon}\int_{A_{4}}\frac{|\xi|^{-\frac{1}{2}}|\xi_{2}|^{\frac{\alpha}{4}}F\prod\limits_{j=1}^{2}F_{j}}
{\langle\sigma_{2}\rangle^{\frac{1}{2}+\epsilon}
   \langle\sigma\rangle^{\frac{1}{2}+\epsilon}}
   d\xi_{1}d\mu_{1}d\tau_{1}d\xi d\mu d\tau\nonumber\\
&&\leq CN^{-\frac{3\alpha}{4}+1+(2\alpha+2)\epsilon}
\|F\|_{L_{\tau\xi\mu}^{2}}\left(\prod_{j=1}^{2}\|F_{j}\|_{L_{\tau\xi\mu}^{2}}\right).
\end{eqnarray*}
 When (\ref{3.018}) is valid,  this case can be proved similarly to case (\ref{3.016}).

\noindent When $|\xi|>6|\xi_{2}|,$ we have that $|\xi|\sim |\xi_{1}|,$ we consider (\ref{3.016})-(\ref{3.018}), respectively.

\noindent When (\ref{3.016}) is valid, since $-\frac{3\alpha}{8}+\frac{1}{2}+4\alpha\epsilon\leq s\leq0,$ we have that
\begin{eqnarray}
   &&\hspace{-1.8cm} K_{2}(\xi_{1},\mu_1,\tau_{1},\xi,\mu,\tau)\leq CN^{2s}\frac{|\xi|\prod\limits_{j=1}^{2}|\xi_{j}|^{-s}}{\langle\sigma\rangle^{\frac{1}{2}-2\epsilon}
   \langle\sigma_{1}\rangle^{\frac{1}{2}+\epsilon}\langle\sigma_{2}\rangle^{\frac{1}{2}+\epsilon}}\nonumber\\&&\leq C N^{2s}\frac{|\xi_{2}|^{-\frac{1}{2}+2\epsilon-s}|\xi_{1}|^{1-s-\frac{\alpha}{2}+2\alpha\epsilon}}{\langle\sigma_{1}\rangle^{\frac{1}{2}+\epsilon}
   \langle\sigma_{2}\rangle^{\frac{1}{2}+\epsilon}}\nonumber\\&&\leq C
   |\xi_{1}|^{-2s-\frac{3\alpha}{4}+1+(2\alpha+2)\epsilon}\frac{|\xi_{2}|^{-\frac{1}{2}}
   |\xi_{1}|^{\frac{\alpha}{4}}}{\langle\sigma_{1}\rangle^{\frac{1}{2}+\epsilon}
   \langle\sigma_{2}\rangle^{\frac{1}{2}+\epsilon}}\nonumber\\&&\leq CN^{-\frac{3\alpha}{4}+1+(2\alpha+2)\epsilon}\frac{|\xi_{2}|^{-\frac{1}{2}}
   |\xi_{1}|^{\frac{\alpha}{4}}}{\langle\sigma_{1}\rangle^{\frac{1}{2}+\epsilon}
   \langle\sigma_{2}\rangle^{\frac{1}{2}+\epsilon}}.
  \label{3.061}
\end{eqnarray}
Combining (\ref{2.011}) with (\ref{3.061}), we have that
\begin{eqnarray}
&&J_{4}\leq CN^{-\frac{3\alpha}{4}+1+(2\alpha+2)\epsilon}\int_{A_{4}}\frac{|\xi_{2}|^{-\frac{1}{2}}|\xi_{1}|^{\frac{\alpha}{4}}F\prod\limits_{j=1}^{2}F_{j}}{\langle\sigma_{1}\rangle^{\frac{1}{2}+\epsilon}
   \langle\sigma_{2}\rangle^{\frac{1}{2}+\epsilon}}d\xi_{1}d\mu_{1}d\tau_{1}d\xi d\mu d\tau\nonumber\\&&\leq CN^{-\frac{3\alpha}{4}+1+(2\alpha+2)\epsilon}\|F\|_{L_{\tau\xi\mu}^{2}}\left(\prod_{j=1}^{2}\|F_{j}\|_{L_{\tau\xi\mu}^{2}}\right).\label{3.062}
\end{eqnarray}
When (\ref{3.017}) is valid, since $\langle \sigma\rangle ^{-\frac{1}{2}+2\epsilon}\langle \sigma_{1}\rangle^{-\frac{1}{2}-\epsilon}\leq \langle \sigma_{1}\rangle ^{-\frac{1}{2}+2\epsilon}\langle \sigma\rangle^{-\frac{1}{2}-\epsilon}$ and  $-\frac{3\alpha}{8}+\frac{1}{2}+4\alpha\epsilon\leq s\leq0,$ we have that
\begin{eqnarray}
   && K_{2}(\xi_{1},\mu_1,\tau_{1},\xi,\mu,\tau)\leq CN^{2s}\frac{|\xi|\prod\limits_{j=1}^{2}|\xi_{j}|^{-s}}{\langle\sigma_{1}\rangle^{\frac{1}{2}-2\epsilon}
   \langle\sigma\rangle^{\frac{1}{2}+\epsilon}\langle\sigma_{2}\rangle^{\frac{1}{2}+\epsilon}}\nonumber\\&&\leq C N^{2s}\frac{|\xi_{2}|^{-\frac{1}{2}+2\epsilon-s}|\xi_{1}|^{1-s-\frac{\alpha}{2}+2\alpha+\epsilon}}{\langle\sigma\rangle^{\frac{1}{2}+\epsilon}
   \langle\sigma_{2}\rangle^{\frac{1}{2}+\epsilon}}\leq CN^{2s} \frac{|\xi_{2}|^{-\frac{1}{2}}
   |\xi_{1}|^{\frac{\alpha}{4}}|\xi_{1}|^{1-2s-\frac{3\alpha}{4}
   +(2\alpha+2)\epsilon}}{\langle\sigma\rangle^{\frac{1}{2}+\epsilon}
   \langle\sigma_{2}\rangle^{\frac{1}{2}+\epsilon}}\nonumber\\&&\leq CN^{-\frac{3\alpha}{4}+1+(2\alpha+2)\epsilon}\frac{|\xi_{2}|^{-\frac{1}{2}}
   |\xi|^{\frac{\alpha}{4}}}{\langle\sigma\rangle^{\frac{1}{2}+\epsilon}
   \langle\sigma_{2}\rangle^{\frac{1}{2}+\epsilon}}.
  \label{3.063}
\end{eqnarray}
Combining (\ref{2.011}) with (\ref{3.063}), we have that
\begin{eqnarray}
&&J_{4}\leq CN^{-\frac{3\alpha}{4}+1+(2\alpha+2)\epsilon}\int_{A_{4}}\frac{|\xi_{2}|^{-\frac{1}{2}}|\xi|^{\frac{\alpha}{4}}F\prod\limits_{j=1}^{2}F_{j}}
{\langle\sigma\rangle^{\frac{1}{2}+\epsilon}
   \langle\sigma_{2}\rangle^{\frac{1}{2}+\epsilon}}d\xi_{1}d\mu_{1}d\tau_{1}d\xi d\mu d\tau\nonumber\\&&\leq CN^{-\frac{3\alpha}{4}+1+(2\alpha+2)\epsilon}\|F\|_{L_{\tau \xi\mu}^{2}}\left(\prod_{j=1}^{2}\|F_{j}\|_{L_{\tau\xi\mu}^{2}}\right).\label{3.064}
\end{eqnarray}
When (\ref{3.018}) is valid, since $\langle \sigma\rangle ^{-\frac{1}{2}+2\epsilon}\langle \sigma_{2}\rangle^{-\frac{1}{2}-\epsilon}\leq \langle \sigma_{2}\rangle ^{-\frac{1}{2}+2\epsilon}\langle \sigma\rangle^{-\frac{1}{2}-\epsilon}$ and  $-\frac{3\alpha}{8}+\frac{1}{2}+4\alpha\epsilon\leq s\leq0,$ we have that
\begin{eqnarray}
  && K_{2}(\xi_{1},\mu_1,\tau_{1},\xi,\mu,\tau)\leq CN^{2s}\frac{|\xi|\prod\limits_{j=1}^{2}|\xi_{j}|^{-s}}{\langle\sigma_{2}\rangle^{\frac{1}{2}-2\epsilon}
   \langle\sigma\rangle^{\frac{1}{2}+\epsilon}\langle\sigma_{1}\rangle^{\frac{1}{2}+\epsilon}}\nonumber\\&&\leq C N^{2s}\frac{|\xi_{2}|^{-\frac{1}{2}+2\epsilon-s}|\xi_{1}|^{1-s-\frac{\alpha}{2}+2\alpha\epsilon}}{\langle\sigma\rangle^{\frac{1}{2}+\epsilon}
   \langle\sigma_{1}\rangle^{\frac{1}{2}+\epsilon}};
  \label{3.065}
\end{eqnarray}
if $-\frac{1}{2}+2\epsilon\leq s\leq0,$ from (\ref{3.065}), we have that
\begin{eqnarray}
  && K_{2}(\xi_{1},\mu_1,\tau_{1},\xi,\mu,\tau)\leq C N^{s-\frac{1}{2}+2\epsilon}\frac{|\xi|^{-\frac{1}{4}+\frac{\alpha}{8}}|\xi_{1}|^{-\frac{1}{4}+\frac{\alpha}{8}}
  |\xi_{1}|^{-\frac{3\alpha}{4}+\frac{3}{2}-s+2\alpha\epsilon}}{\langle\sigma\rangle^{\frac{1}{2}+\epsilon}
   \langle\sigma_{1}\rangle^{\frac{1}{2}+\epsilon}}\nonumber\\&&\leq CN^{-\frac{3\alpha}{4}+1+(2\alpha+2)\epsilon}\frac{|\xi|^{-\frac{1}{4}+\frac{\alpha}{8}}|\xi_{1}|^{-\frac{1}{4}+\frac{\alpha}{8}}
  }{\langle\sigma\rangle^{\frac{1}{2}+\epsilon}
   \langle\sigma_{1}\rangle^{\frac{1}{2}+\epsilon}};
  \label{3.066}
\end{eqnarray}
if $-\frac{3\alpha}{8}+\frac{1}{2}+4\alpha\epsilon\leq s<-\frac{1}{2}+2\epsilon,$ from (\ref{3.065}), we have that
\begin{eqnarray}
  && K_{2}(\xi_{1},\mu_1,\tau_{1},\xi,\mu,\tau)\leq C N^{2s}\frac{|\xi|^{-\frac{1}{4}+\frac{\alpha}{8}}|\xi_{1}|^{-\frac{1}{4}+\frac{\alpha}{8}}
  |\xi_{1}|^{-\frac{3\alpha}{4}+1-2s+2\alpha\epsilon}}{\langle\sigma\rangle^{\frac{1}{2}+\epsilon}
   \langle\sigma_{1}\rangle^{\frac{1}{2}+\epsilon}}\nonumber\\&&\leq CN^{-\frac{3\alpha}{4}+1+(2\alpha+2)\epsilon}\frac{|\xi|^{-\frac{1}{4}+\frac{\alpha}{8}}|\xi_{1}|^{-\frac{1}{4}+\frac{\alpha}{8}}
  }{\langle\sigma\rangle^{\frac{1}{2}+\epsilon}
   \langle\sigma_{1}\rangle^{\frac{1}{2}+\epsilon}}.
  \label{3.067}
\end{eqnarray}
Combining (\ref{2.09}) with (\ref{3.066})-(\ref{3.067}), we have that
\begin{eqnarray}
&&J_{4}\leq CN^{-\frac{3\alpha}{4}+1+(2\alpha+2)\epsilon}\int_{A_{4}}\frac{|\xi_{1}|^{-\frac{1}{2}}|\xi|^{\frac{\alpha}{4}}F\prod\limits_{j=1}^{2}F_{j}}
{\langle\sigma\rangle^{\frac{1}{2}+\epsilon}
   \langle\sigma_{1}\rangle^{\frac{1}{2}+\epsilon}}d\xi_{1}d\mu_{1}d\tau_{1}d\xi d\mu d\tau\nonumber\\&&\leq CN^{-\frac{3\alpha}{4}+1+(2\alpha+2)\epsilon}\|F\|_{L_{\tau\xi\mu}^{2}}\left(\prod_{j=1}^{2}\|F_{j}\|_{L_{\tau\xi\mu}^{2}}\right).\label{3.068}
\end{eqnarray}

This completes the proof of Lemma 3.2.

\begin{Lemma}\label{Lem3.3}
Let $s\geq -\frac{3\alpha}{8}+\frac{1}{4}+4\alpha\epsilon$ and $u_{j}\in X_{\frac{1}{2}+\epsilon}^{s,0}(j=1,2)$.
Then, we have that
\begin{eqnarray}
&&\|\partial_{x}I(u_{1}u_{2})\|_{X_{-\frac{1}{2}+2\epsilon}^{0,0}}\leq C
\prod_{j=1}^{2}\|Iu_{j}\|_{X_{\frac{1}{2}+\epsilon}^{0,0}}.\label{3.069}
\end{eqnarray}\end{Lemma}\noindent{\bf Proof.} To prove (\ref{3.069}),  by duality, it suffices to  prove that
\begin{eqnarray}
&&\left|\int_{\SR^{3}}\bar{u}\partial_{x}I(u_{1}u_{2})dxdydt\right|\leq
C\|u\|_{X_{\frac{1}{2}-2\epsilon}^{0,0}}\left(\prod_{j=1}^{2}
\|Iu_{j}\|_{X_{\frac{1}{2}+\epsilon}^{0,0}}\right).\label{3.070}
\end{eqnarray}
for $u\in X_{\frac{1}{2}-2\epsilon}^{0,0}.$
Let
\begin{eqnarray}
&&\xi=\xi_{1}+\xi_{2},\mu=\mu_{1}+\mu_{2},\tau=\tau_{1}+\tau_{2},\nonumber\\
&&F(\xi,\mu,\tau)=
\langle \sigma\rangle^{\frac{1}{2}-2\epsilon}\mathscr{F}u(\xi,\mu,\tau),\nonumber\\&&
F_{j}(\xi_{j},\mu_{j},\tau_{j})=M(\xi_{j})
\langle \sigma_{j}\rangle^{\frac{1}{2}+\epsilon}
\mathscr{F}u_{j}(\xi_{j},\mu,\tau_{j})(j=1,2).\label{3.071}
\end{eqnarray}
To obtain (\ref{3.070}), from (\ref{3.071}), it suffices to prove that
\begin{eqnarray}
&&\int_{\SR^{6}}\frac{|\xi|M(\xi)
F(\xi,\mu,\tau)\prod\limits_{j=1}^{2}F_{j}(\xi_{j},\mu_{j},\tau_{j})}{\langle\sigma\rangle^{\frac{1}{2}-2\epsilon}
\prod\limits_{j=1}^{2}M(\xi_{j})\langle\sigma_{j}\rangle^{\frac{1}{2}+\epsilon}}
d\xi_{1}d\mu_{1}d\tau_{1}d\xi d\mu d\tau\nonumber\\&&\leq C
\|F\|_{L_{\tau\xi\mu}^{2}}\left(\prod_{j=1}^{2}\|F_{j}\|_{L_{\tau\xi\mu}^{2}}\right).\label{3.072}
\end{eqnarray}
From (2.4) of \cite{IMEJDE}, we have that
\begin{eqnarray}
\frac{M(\xi)}{\prod\limits_{j=1}^{2}M(\xi_{j})}\leq C\frac{\langle\xi\rangle^{s}}
{\prod\limits_{j=1}^{2}\langle\xi_{j}\rangle^{s}}\label{3.073}.
\end{eqnarray}
By using (\ref{3.073}), we have that the left hand side of (\ref{3.072})  can be bounded by
\begin{eqnarray}
&&\int_{\SR^{6}}\frac{|\xi|\langle\xi\rangle^{s}
F(\xi,\mu,\tau)\prod\limits_{j=1}^{2}F_{j}(\xi_{j},\mu_{j},\tau_{j})}{\langle\sigma_{j}\rangle^{\frac{1}{2}-2\epsilon}
\prod\limits_{j=1}^{2}\langle\xi_{j}\rangle^{s}\langle\sigma_{j}\rangle^{\frac{1}{2}+\epsilon}}
d\xi_{1}d\mu_{1}d\tau_{1}d\xi d\mu d\tau.\label{3.074}
\end{eqnarray}
By using (\ref{3.01}),  we have that (\ref{3.074}) can be bounded by
$ C
\|F\|_{L_{\tau\xi\mu}^{2}}\left(\prod\limits_{j=1}^{2}\|F_{j}\|_{L_{\tau\xi\mu}^{2}}\right).$

This completes the proof of Lemma 3.3.

\bigskip
\bigskip
\noindent {\large\bf 4. Proof of Theorem  1.1}

\setcounter{equation}{0}

 \setcounter{Theorem}{0}

\setcounter{Lemma}{0}

\setcounter{section}{4}

In this section, combining Lemmas 2.2, 3.1 with the fixed point theorem, we present the proof of Theorem 1.1.

\noindent We define
\begin{eqnarray}
&&\Phi_{1}(u):=\psi(t)W(t)u_{0}+\frac{1}{2}\psi\left(\frac{t}{T}\right)
\int_{0}^{t}W(t-\tau)\partial_{x}(u^{2})d\tau,\label{4.01}\\
&&B_{1}(0,2C\|u_{0}\|_{H^{s_{1},s_{2}}}):=\left\{u:\|u\|_{X_{\frac{1}{2}
+\epsilon}^{s_{1},s_{2}}}\leq 2C\|u_{0}\|_{H^{s_{1},s_{2}}}\right\}.\label{4.02}
\end{eqnarray}
Here $\psi(t)$ be defined as in line 2 from bottom of page 5.
Combining Lemmas 2.2, 3.1 with (\ref{4.01}), (\ref{4.02}),  we have that
\begin{eqnarray}
&&\left\|\Phi_{1}(u)\right\|_{X_{\frac{1}{2}+\epsilon}^{s_{1},s_{2}}}\leq \left\|\psi(t)W(t)u_{0}\right\|_{X_{\frac{1}{2}+\epsilon}^{s_{1},s_{2}}}
+\left\|\frac{1}{2}\psi\left(\frac{t}{T}\right)\int_{0}^{t}W(t-\tau)\partial_{x}(u^{2})d\tau\right\|_{X_{\frac{1}{2}+\epsilon}^{s_{1},s_{2}}}\nonumber\\
&&\leq C\|u_{0}\|_{H^{s_{1},s_{2}}}+CT^{\epsilon}\left\|\partial_{x}(u^{2})\right\|_{X_{-\frac{1}{2}+2\epsilon}^{s_{1},s_{2}}}\nonumber\\
&&\leq C\|u_{0}\|_{H^{s_{1},s_{2}}}+CT^{\epsilon}\left\|u\right\|_{X_{\frac{1}{2}+\epsilon}^{s_{1},s_{2}}}^{2}\nonumber\\
&&\leq C\|u_{0}\|_{H^{s_{1},s_{2}}}+4C^{3}T^{\epsilon}\left\|u_{0}\right\|_{H^{s_{1},s_{2}}}^{2}.\label{4.03}
\end{eqnarray}
We choose $T\in (0,1)$ such that
\begin{eqnarray}
T^{\epsilon}=\left[16C^{2}(\|u_{0}\|_{H^{s_{1},s_{2}}}+1)\right]^{-1}.\label{4.04}
\end{eqnarray}
Combining (\ref{4.03}) with (\ref{4.04}),  we have that
\begin{eqnarray}
&&\left\|\Phi_{1}(u)\right\|_{X_{\frac{1}{2}+\epsilon}^{s_{1},s_{2}}}
\leq C\|u_{0}\|_{H^{s_{1},s_{2}}}+C\left\|u_{0}\right\|_{H^{s_{1},s_{2}}}=2C\left\|u_{0}\right\|_{H^{s_{1},s_{2}}}.\label{4.05}
\end{eqnarray}
Thus, $\Phi_{1}$ maps $B_{1}(0,2C\|u_{0}\|_{H^{s_{1},s_{2}}})$ into $B_{1}(0,2C\|u_{0}\|_{H^{s_{1},s_{2}}})$.
By using Lemmas 2.2, 3.1, (\ref{4.04})-(\ref{4.05}), we have that
\begin{eqnarray}
&&\left\|\Phi_{1}(u_{1})-\Phi_{1}(u_{2})\right\|_{X_{\frac{1}{2}+\epsilon}^{s_{1},s_{2}}}
\leq C\left\|\frac{1}{2}\psi\left(\frac{t}{\tau}\right)
\int_{0}^{t}W(t-\tau)\partial_{x}(u_{1}^{2}-u_{2}^{2})d\tau\right\|_{X_{\frac{1}{2}+\epsilon}^{s_{1},s_{2}}}\nonumber\\
&&\leq CT^{\epsilon}\left\|u_{1}-u_{2}\right\|_{X_{\frac{1}{2}+\epsilon}^{s_{1},s_{2}}}
\left[\left\|u_{1}\right\|_{X_{\frac{1}{2}+\epsilon}^{s_{1},s_{2}}}+\left\|u_{2}\right\|_{X_{\frac{1}{2}+\epsilon}^{s_{1},s_{2}}}\right]\nonumber\\
&&\leq 4C^{2}T^{\epsilon}\left\|u_{0}\right\|_{H^{s_{1},s_{2}}}\left\|u_{1}-u_{2}\right\|_{X_{\frac{1}{2}+\epsilon}^{s_{1},s_{2}}}\leq \frac{1}{2}\left\|u_{1}-u_{2}\right\|_{X_{\frac{1}{2}+\epsilon}^{s_{1},s_{2}}}.\label{4.06}
\end{eqnarray}
Thus, $\Phi_{1}$ is a contraction in the closed ball $B_{1}(0,2C\|u_{0}\|_{H^{s_{1},s_{2}}})$.
 Consequently, $u$ is the fixed point of $\Phi$ in the closed ball
$B_{1}(0,2C\|u_{0}\|_{H^{s_{1},s_{2}}})$. Then $v:=u|_{[0,T]}\in X_{\frac{1}{2}+\epsilon}^{s_{1},s_{2}}([0,T])$
is a solution in the interval $[0,T]$ of the Cauchy problem for (\ref{1.01}) with the initial data $u_{0}$.
 For the facts that  uniqueness of the solution and
the solution to the Cauchy problem for (\ref{1.01}) is continuous
with respect to the initial data, we refer the readers  to Theorems II, III of \cite{IMS}.

This ends the proof of Theorem 1.1.
\bigskip
\bigskip

\noindent {\large\bf 5. Proof of Theorem  1.2}

\setcounter{equation}{0}

 \setcounter{Theorem}{0}

\setcounter{Lemma}{0}

\setcounter{section}{5}
In this section, we give the proof of Theorem 1.2. We present the proof of  Lemma 5.1 before giving the proof of Theorem 1.2.

\begin{Lemma}\label{Lemma5.1}
Let $s_{1}>\frac{1}{4}-\frac{3}{8}\alpha$ and $R:=\frac{1}{8(C+1)^{3}}$, where  $C$  is the largest of those constants which appear in (\ref{2.05})-(\ref{2.06}),  (\ref{3.069}).
  Then, the Cauchy problem for (\ref{1.01}) locally well-posed for data satisfying $I_{N}u_{0}\in L^{2}(\R^{2})$ with
  \begin{eqnarray}
 \left\|I_{N}u_{0}\right\|_{L^{2}}\leq R.\label{5.01}
  \end{eqnarray}
Moreover, the solution to the Cauchy problem for (\ref{1.01}) exists on a time interval $[0,1]$.
\end{Lemma}
\noindent{\bf Proof.} We define $v:=I_{N}u$. If $u$ is the solution to the Cauchy problem for (\ref{1.01}), then
$v$ satisfies the following equation
\begin{eqnarray}
v_{t}-|D_{x}|^{\alpha}u_{x}+\partial_{x}^{-1}
\partial_{y}^{2}v+\frac{1}{2}I_{N}\partial_{x}(I_{N}^{-1}v)^{2}=0.\label{5.02}
\end{eqnarray}
Then $v$ is formally equivalent to the following integral equation
\begin{eqnarray}
v=W(t)I_{N}u_{0}+\frac{1}{2}\int_{0}^{t}W(t-\tau)I_{N}\partial_{x}(I_{N}^{-1}v)^{2}.\label{5.03}
\end{eqnarray}
We define
\begin{eqnarray}
\Phi_{2}(v)=\psi(t)W(t)I_{N}u_{0}+\frac{1}{2}\psi(t)
\int_{0}^{t}W(t-\tau)I_{N}\partial_{x}(I_{N}^{-1}v)^{2}.\label{5.04}
\end{eqnarray}
By using Lemmas 2.2,  3.3, we have that
\begin{eqnarray}
&&\left\|\Phi_{2}(v)\right\|_{X_{\frac{1}{2}+\epsilon}^{0,0}}\leq
 \left\|\psi(t)W(t)I_{N}u_{0}\right\|_{X_{\frac{1}{2}+\epsilon}^{0,0}}
+C\left\|\psi(t)\int_{0}^{t}W(t-\tau)I_{N}\partial_{x}
(I_{N}^{-1}v)^{2}\right\|_{X_{\frac{1}{2}+\epsilon}^{0,0}}\nonumber\\
&&\leq C\left\|I_{N}u_{0}\right\|_{L^{2}}+C\left\|I_{N}\partial_{x}
(I_{N}^{-1}v)^{2}\right\|_{X_{-\frac{1}{2}+2\epsilon}^{0,0}}\nonumber\\
&&\leq C\left\|I_{N}u_{0}\right\|_{L^{2}}+C\left\|I_{N}\partial_{x}
(I_{N}^{-1}v)^{2}\right\|_{X_{-\frac{1}{2}+2\epsilon}^{0,0}}\nonumber\\
&&\leq C\left\|I_{N}u_{0}\right\|_{L^{2}}+C\|v\|_{X_{\frac{1}{2}+\epsilon}^{0,0}}^{2}\nonumber\\
&&\leq CR+C\|v\|_{X_{\frac{1}{2}+\epsilon}^{0,0}}^{2}.\label{5.05}
\end{eqnarray}
We define
\begin{eqnarray}
B_{2}(0,2CR):=\left\{v:\|v\|_{X_{\frac{1}{2}+\epsilon}^{0,0}}\leq
2CR\right\}.\label{5.06}
\end{eqnarray}
Combining   (\ref{5.05})-(\ref{5.06}) with the definition of  $R$,  we have that
\begin{eqnarray}
&&\left\|\Phi_{2}(v)\right\|_{X_{\frac{1}{2}+\epsilon}^{0,0}}
\leq CR+4C^{3}R^{2}
=2CR.\label{5.07}
\end{eqnarray}
Thus, $\Phi_{2}$ maps $B_{2}(0,2CR)$ into
$B_{2}(0,2CR)$. We define
\begin{eqnarray}
v_{j}=I_{N}u_{j}(j=1,2),w_{1}=I_{N}^{-1}v_{1}-I_{N}^{-1}v_{2},
w_{2}=I_{N}^{-1}v_{1}+I_{N}^{-1}v_{2}.\label{5.08}
\end{eqnarray}
By using Lemmas 2.2, 3.1, 3.2,  (\ref{5.05})-(\ref{5.06}) and the definition of $R$, we have that
\begin{eqnarray}
&&\left\|\Phi_{2}(v_{1})-\Phi_{2}(v_{2})\right\|_{X_{\frac{1}{2}+\epsilon}^{0,0}}
\leq C\left\|\psi(t)
\int_{0}^{t}W(t-\tau)\partial_{x}I_{N}\left[(I_{N}^{-1}v_{1})^{2}-(I_{N}^{-1}v_{2})^{2}\right]
d\tau\right\|_{X_{\frac{1}{2}+\epsilon}^{0,0}}\nonumber\\
&&\leq C\left\|\partial_{x}I_{N}(w_{1}w_{2})\right\|_{X_{-\frac{1}{2}+2\epsilon}^{0,0}}\nonumber\\
&&\leq C\|v_{1}-v_{2}\|_{X_{\frac{1}{2}+\epsilon}^{0,0}}
\left[\|v_{1}\|_{X_{\frac{1}{2}+\epsilon}^{0,0}}+\|v_{2}\|_{X_{\frac{1}{2}+\epsilon}^{0,0}}\right]\nonumber\\
&&\leq 4C^{2}R
\|v_{1}-v_{2}\|_{X_{\frac{1}{2}+\epsilon}^{0,0}}\leq \frac{1}{2}\|v_{1}-v_{2}\|_{X_{\frac{1}{2}+\epsilon}^{0,0}}
.\label{5.09}
\end{eqnarray}
Thus, $\Phi_{2}$ is a contraction in the closed ball $B_{2}(0,2CR)$.
 Consequently, $u$ is the fixed point of $\Phi_{2}$ in the closed ball
$B_{2}(0,2CR)$. Then $v:=u|_{[0,1]}\in X_{\frac{1}{2}+\epsilon}^{0,0}([0,1])$
is a solution in the interval $[0,1]$ of the Cauchy problem for (\ref{5.03}) with the initial data $I_{N}u_{0}$.
For the uniqueness of the solution, we  refer the readers to Theorem II of \cite{IMS}.
 For the fact that
the solution to the Cauchy problem for (\ref{5.03}) is continuous
with respect to the initial data, we refer the readers to Theorem III of  \cite{IMS}. Since the phase function $\phi(\xi,\mu)$
is singular at $\xi=0$, to  define  the  derivative  of  $W(t)u_{0}$,
the requirement $|\xi|^{-1}\mathscr{F}_{xy}u_{0}(\xi,\mu)\in \mathscr{S}^{'}(\R^{2})$  is necessary.

This ends the proof of Lemma 5.1.

Inspired by  \cite{ILM-CPAA}, we use Lemmas 2.7, 3.2, 5.1 to prove Theorem 1.2.

For $\lambda>0$, we define
\begin{eqnarray}
u_{\lambda}(x,y,t)=\lambda^{\alpha}u
\left(\lambda x,\lambda^{\frac{\alpha}{2}+1}y,\lambda^{\alpha+1} t\right),
 u_{0\lambda}(x,y)=\lambda^{\alpha}
u_{0}\left(\lambda x,\lambda^{\frac{\alpha}{2}+1}y\right).\label{5.010}
\end{eqnarray}
Thus, $u_{\lambda}(x,y,t)\in X_{\frac{1}{2}+\epsilon}^{s_{1},0}([0,\frac{T}{\lambda^{\alpha+1}}])$ is the solution to
\begin{eqnarray}
&&\partial_{t}u_{\lambda}-|D_{x}|^{\alpha}\partial_{x}u_{\lambda}
+
\partial_{x}^{-1}\partial_{y}^{2}u_{\lambda}+u_{\lambda}\partial_{x}u_{\lambda}=0,\label{5.011}\\
&&u_{\lambda}(x,y,0)=u_{0\lambda}(x,y),\label{5.012}
\end{eqnarray}
if and only if $u(x,y,t)\in X_{\frac{1}{2}+\epsilon}^{s,0}([0,T])$ is the solution to the
Cauchy problem for (\ref{1.01}) in $[0,T]$ with the initial data $u_{0}.$
By using a direct computation, for $\lambda \in (0,1),$  we have that
\begin{eqnarray}
\|I_{N}u_{0\lambda}\|_{L^{2}}\leq CN^{-s}\lambda^{\frac{3\alpha}{4}-1+s}\|u_{0}\|_{H^{s,0}}.\label{5.013}
\end{eqnarray}
For $u_{0}\neq0$ and $u_{0}\in H^{s,0}(\R^{2})$, we choose $\lambda,N$ such that
\begin{eqnarray}
\|I_{N}u_{0\lambda}\|_{L^{2}}\leq CN^{-s}\lambda^{\frac{3\alpha}{4}-1+s}
\|u_{0}\|_{H^{s,0}}:=\frac{R}{4}.\label{5.014}
\end{eqnarray}
Then there exist $w_{3}$ which satisfies that  $\|w_{3}\|_{X_{\frac{1}{2}+\epsilon}^{s,0}}\leq  2CR$ such that $v:=w_{3}\mid_{[0,1]}$
is a solution to the Cauchy problem for (\ref{5.02}) with $u_{0\lambda}$.
Multiplying (\ref{5.02}) by $2I_{N}u_{\lambda}$ and integrating with respect to $x,y$
and integrating by parts with respect to $x$ yield
\begin{eqnarray}
\frac{d}{dt}\int_{\SR^{2}}(I_{N}u_{\lambda})^{2}dxdy+\int_{\SR^{2}}I_{N}u_{\lambda}\partial_{x}I_{N}
\left[(u_{\lambda})^{2}\right]dxdy=0.\label{5.015}
\end{eqnarray}
Combining
\begin{eqnarray*}
\int_{\SR^{2}}I_{N}u_{\lambda}\partial_{x}\left[(I_{N}u_{\lambda})^{2}\right]dxdy=0
\end{eqnarray*}
with (\ref{5.015}),
we have that
\begin{eqnarray}
\frac{d}{dt}\int_{\SR^{2}}(I_{N}u_{\lambda})^{2}dxdy=-
\int_{\SR^{2}}I_{N}u_{\lambda}\partial_{x}\left[I_{N}\left(u_{\lambda}^{2}\right)-(I_{N}u_{\lambda})^{2}\right]dxdy.\label{5.016}
\end{eqnarray}
From (\ref{5.016}) and Lemma 2.7,  we have that
\begin{eqnarray}
&&\int_{\SR^{2}}(I_{N}u_{\lambda}(x,y,1))^{2}dxdy-\int_{\SR^{2}}(I_{N}u_{0\lambda})^{2}dxdy\nonumber\\&&
=-\int_{0}^{1}\int_{\SR^{2}}I_{N}u_{\lambda}\partial_{x}\left[I_{N}\left((u_{\lambda})^{2}
\right)-(I_{N}u_{\lambda})^{2}\right]dxdydt\nonumber\\
&&=-\int_{\SR}\int_{\SR^{2}}\left(\chi_{[0,1]}(t)I_{N}u_{\lambda}\right)
\left(\chi_{[0,1]}(t)\partial_{x}\left[I_{N}\left((u_{\lambda})^{2}\right)-(I_{N}u_{\lambda})^{2}\right]\right)dxdydt
\nonumber\\
&&\leq C\left\|\chi_{[0,1]}(t)I_{N}u_{\lambda}\right\|_{X_{\frac{1}{2}-\epsilon}^{0,0}}
\left\|\chi_{[0,1]}(t)\partial_{x}\left[I_{N}\left((u_{\lambda})^{2}\right)-(I_{N}u_{\lambda})^{2}\right]
\right\|_{X_{-\frac{1}{2}+\epsilon}^{0,0}}\nonumber\\
&&\leq C\left\|I_{N}u_{\lambda}\right\|_{X_{\frac{1-\epsilon}{2}}^{0,0}}
\left\|\partial_{x}\left[I_{N}\left((u_{\lambda})^{2}\right)-(I_{N}u_{\lambda})^{2}\right]\right\|_{X_{-\frac{1}{2}+2\epsilon}^{0,0}}\nonumber\\
&&\leq CN^{-\frac{3\alpha}{4}+1+(2\alpha+2)\epsilon}\|I_{N}u_{\lambda}\|_{X_{\frac{1}{2}+\epsilon}^{0,0}}^{3}.\label{5.017}
\end{eqnarray}
Combining  (\ref{5.014}) with  (\ref{5.017}), we have that
\begin{eqnarray}
&&\int_{\SR^{2}}(I_{N}u_{\lambda}(x,y,1))^{2}dxdy\leq \frac{R^{2}}{4}+ CN^{-\frac{3\alpha}{4}+1+(2\alpha+2)\epsilon}
\|I_{N}u_{\lambda}\|_{X_{\frac{1}{2}+\epsilon}^{0,0}}^{3}\nonumber\\
&&\leq \frac{R^{2}}{4}+8C^{4}N^{-\frac{3\alpha}{4}+1+(2\alpha+2)\epsilon}R^{3}.\label{5.018}
\end{eqnarray}
Thus, if we take $N$ sufficiently large such that  $8C^{4}N^{-\frac{3\alpha}{4}+1+(2\alpha+2)\epsilon}R^{3}\leq \frac{3}{4}R^{2},$
then
\begin{eqnarray}
\left[\int_{\SR^{2}}(I_{N}u_{\lambda}(x,y,1))^{2}dxdy\right]^{\frac12}\leq R.\label{5.019}
\end{eqnarray}
We consider $I_{N}u_{\lambda}(x,y,1)$ as the initial data and repeat the above argument, we obtain that (\ref{5.011})-(\ref{5.012}) possess
a solution in $\R^{2}\times [1,2]$. In this way, we can extend the solution to (\ref{5.011})-(\ref{5.012}) to the time interval
$[0,2].$ The above argument can be repeated $L$ steps, where $L$ is the maximal positive integer  such that
\begin{eqnarray}
8C^{4}N^{-\frac{3\alpha}{4}+1+(2\alpha+2)\epsilon}R^{3}L\leq \frac{3}{4}R^{2}.\label{5.020}
\end{eqnarray}
More precisely, the solution to (\ref{5.011})-(\ref{5.012}) can be extended to the time interval $[0,L]$. Thus, we can prove that (\ref{5.011})-(\ref{5.012})
are globally well-posed in $[0,\frac{T}{\lambda^{\alpha+1}}]$ if  we can choose a number $N$ such that
\begin{eqnarray}
L\geq \frac{T}{\lambda^{\alpha+1}}.\label{5.021}
\end{eqnarray}
From (\ref{5.020}), we know that
\begin{eqnarray}
L\sim N^{-\frac{3\alpha}{4}+1+(2\alpha+2)\epsilon}.\label{5.022}
\end{eqnarray}
We know  that (\ref{5.021}) is valid provided that the following inequality is valid
\begin{eqnarray}
CN^{\frac{3\alpha}{4}-1-(2\alpha+2)\epsilon}\geq \frac{T}{\lambda^{\alpha+1}}\sim CTN^{\frac{-(\alpha+1)s}{\frac{3}{4}\alpha-1+s}}.\label{5.023}
\end{eqnarray}
In fact,  (\ref{5.023}) is valid if
\begin{eqnarray}
N^{\frac{3\alpha}{4}-1}> N^{\frac{-(\alpha+1)s}{\frac{3}{4}\alpha-1+s}}\label{5.024}
\end{eqnarray}
which is equivalent to $-\frac{(3\alpha-4)^{2}}{28\alpha}<s_{1}<0.$

This completes the proof of Theorem 1.2.

\bigskip
\bigskip

\leftline{\large \bf Acknowledgments}

\bigskip

\noindent

 This work is supported by the Natural Science Foundation of China
 under grant numbers 11401180 and 11171116 and 11471330. The first author is also
 supported by   the Young core Teachers program of Henan Normal University and  15A110033.

\end{document}